\documentclass[a4paper]{amsart}
\usepackage{amsfonts}
\usepackage{amssymb}
\usepackage[utf8]{inputenc}
\usepackage{amsmath}
\usepackage{pdflscape}
\usepackage{graphicx}%

\providecommand{\U}[1]{\protect\rule{.1in}{.1in}}
\usepackage[colorlinks=true, linkcolor=red, citecolor=blue]{hyperref}
\usepackage[]{epsfig}
\usepackage[]{pstricks}
\usepackage{tikz}
\newtheorem{theorem}{Theorem}[section]
\newtheorem{lemma}{Lemma}[section]

\newtheorem{remark*}{Remark}

\newtheorem{remarks*}{Remarks}

\newtheorem{claim}{Claim}
\numberwithin{equation}{section}
\newcommand{\R}{\mathbb R}

\newcommand{\N}{{\mathbb N}}

\newcommand{\ft}{{\mathcal{F}}}

\newcommand{\Sch}{{\mathcal{S}}}
\newcommand{\supp}{{\mbox{supp}}}

\newcommand{\px}{\partial_x}
\newcommand{\pt}{\partial_t}

\def\norm#1{\|#1\|}
\def\bra#1{\langle#1\rangle}
\def\wt#1{\widetilde{#1}}
\def\wh#1{\widehat{#1}}
\def\set#1{\{#1\}}

\providecommand{\norm}[1]{\left\lVert#1\right\rVert}
\numberwithin{equation}{section}

\begin{document}

	\pagenumbering{arabic}	
\title[4NLS on graphs]{Forcing operators on star graphs applied for the cubic fourth order Schr\"odinger equation}
\author[Capistrano--Filho]{Roberto A. Capistrano--Filho*}\thanks{*Corresponding author: roberto.capistranofilho@ufpe.br ; capistranofilho@dmat.ufpe.br}
\address{Departamento de Matem\'atica,  Universidade Federal de Pernambuco (UFPE), 50740-545, Recife (PE), Brazil.}
\email{roberto.capistranofilho@ufpe.br ; capistranofilho@dmat.ufpe.br}
\author[Cavalcante]{M\'arcio Cavalcante}
\address{Instituto de Matem\'{a}tica, Universidade Federal de Alagoas (UFAL), Macei\'o (AL), Brazil}
\email{marcio.melo@im.ufal.br}
\author[Gallego]{Fernando A. Gallego}
\address{Departamento de Matematicas y Estad\'istica, Universidad Nacional de Colombia (UNAL), Cra 27 No. 64-60, 170003, Manizales, Colombia}
\email{fagallegor@unal.edu.co}
%\thanks{}
\subjclass[2010]{35R02, 35Q55, 35C15, 81Q35, 35G30} 
\keywords{Unbounded star graphs structure, Fourth order Schr\"odinger equation, Boundary forcing operator approach, Boundary conditions}
%\date{Version 2019-09-15}

\begin{abstract}
In a recent article \textit{``Lower regularity solutions of the biharmonic Schrödinger equation in a quarter plane'', to appear on Pacific Journal of Mathematics \cite{CaCaGa}}, the authors gave a starting point of the study on a series of problems concerning the initial boundary value problem and control theory of Biharmonic NLS in some non-standard domains. In this direction, this article deals to present answers for some questions left in \cite{CaCaGa} concerning the study of the cubic fourth order Schr\"odinger equation in a star graph structure $\mathcal{G}$. Precisely, consider $\mathcal{G}$ composed by $N$ edges parameterized by half-lines $(0,+\infty)$ attached with a common vertex $\nu$. With this structure the manuscript proposes to study the well-posedness of a dispersive model on star graphs with three appropriated vertex conditions by using the \textit{boundary forcing operator approach}. More precisely, we give positive answer for the Cauchy problem in low regularity Sobolev spaces. We have noted that this approach seems very efficient, since this allows to use the tools of Harmonic Analysis,  for instance, the Fourier restriction method, introduced by Bourgain, while for the other known standard methods to solve partial differential partial equations on star graphs are more complicated to capture the dispersive smoothing effect in low regularity.  The arguments presented in this work have prospects to be applied for other nonlinear dispersive equations in the context of star graphs with unbounded edges.
\end{abstract}

\maketitle

\section{Introduction}
\subsection{Quantum and metric graphs}

In mathematics and physics, a quantum graph is a linear network-shaped structure of vertices connected on edges (i.e., a graph), where a differential (or pseudo-differential) equation is posed on each edge, while in the case of each edge is equipped with a natural metric  the graph is denoted as a metric graph. An example would be a power network consisting of power lines (edges) connected at transformer stations (vertices); the differential equations would be  then the voltage along each of the line and the boundary conditions for each edge equipped at the adjacent vertices ensuring that the current added over all edges adds to zero at each vertex.

%\textcolor{red}{In other words, we can express a quantum graph as a metric graph or a network-shaped structure of vertices connected by edges, with a linear Hamiltonian operator (such as a Schr\"odinger-like operator) suitably defined on functions that are supported on the edges. It arises as a simplified models in for wave propagation, for instance, in a  quasi one-dimensional  (e.g. meso - or nanoscale) system that looks like a thin  neighborhood of a  graph. 	} 		

Quantum graphs were first studied by Linus Pauling as models of free electrons in organic molecules in the 1930s. They also appear in a variety of mathematical contexts, e.g. as model systems in quantum chaos, in the study of waveguides, in photonic crystals and in Anderson localization - is the absence of diffusion of waves in a disordered medium -, or as limit on shrinking thin wires. Quantum graphs have become prominent models in mesoscopic physics used to obtain a theoretical understanding of nanotechnology. Another, more simple notion of quantum graphs was introduced by Freedman \textit{et al.} in \cite{Freedman}.

Aside from actually solving the differential equations posed on a quantum graph for purposes of concrete applications, typical questions that arise are those of well-posedness, controllability (what inputs have to be provided to bring the system into a desired state, for example providing sufficient power to all houses on a power network) and identifiability (how and where one has to measure something to obtain a complete picture of the state of the system, for example measuring the pressure of a water pipe network to determine whether or not there is a leaking pipe).

\subsection{Nonlinear dispersive models on star graphs}
In the last years, the study of nonlinear dispersive models in a metric graph has attracted a lot of attention of mathematicians, physicists, chemists and engineers,  see for details \cite{BK, BlaExn08, BurCas01, K, Mug15} and references therein. In particular, the framework prototype  (graph-geometry)  for description of these phenomena have been a {\it star graph} $\mathcal G$, namely, on metric graphs with $N$  half-lines of the form $(0, +\infty)$  connecting at a common vertex $\nu=0$, together with a nonlinear equation suitably defined on the edges such as the nonlinear Schr\"odinger equation (see Adami {\it{et al.}} \cite{AdaNoj14, AdaNoj15} and Angulo and Goloshchapova \cite{AngGol17a, AngGol17b}). We note that with the introduction of nonlinearities in the dispersive models, the network provides a nice field, where one can look for interesting soliton propagation and nonlinear dynamics in general.   A central point that makes this analysis a delicate problem is the presence of a vertex where the underlying one-dimensional star graph should bifurcate (or multi-bifurcate in a general metric graph). 

%\textcolor{blue}{In others contexts, for example, in the context of soliton transport in networks and branched structures (see \cite{SBM, SMSSK}) since wave dynamics in networks can be modeled by nonlinear evolution equations suitably defined on the edges. Soliton and other nonlinear waves in branched systems appear in different system of condensed matter, Josephson junction networks, polymers, optics, neuroscience, DNA, blood pressure waves in large arteries or in shallow water equation to describe a fluid network (see e.g. \cite{AdaNoj13a, Berkolaiko, BK, BeK, BurCas01, CM, Fid15, K,  Mug15,Mehmeti, Noj14}).} 

%To address these issues, in general the problem is difficult to tackle because both  equations of motion and the geometry are complex.

%We note that not branching angles but the topology of bifurcation is essential. Indeed,  a soliton-profile coming into the vertex along one of the bonds shows a complicated motion around the vertex such as reflection and emergence of the radiation there, moreover, in particular one cannot see easily how energy travels across the network. Therefore,  the study of the dynamic for non-linear evolution models becomes a challenge and it will depend  heavily on the conditions on the vertex (or vertices) to have a fruitful description of the system. 

Looking at other nonlinear dispersive systems on graphs structure, we have some interesting results. For example, related with well-posedness theory, the second author in \cite{Cav}, studied the local well-posedness for the Cauchy problem associated to Korteweg-de Vries equation in a metric star graph with three semi-infinite edges given by one negative half-line and two positives half-lines attached to a common vertex $\nu=0$ (the $\mathcal Y$-junction framework).  Another nonlinear dispersive equation, the Benjamin--Bona--Mahony (BBM) equation, is treated in \cite{bona,Mugnolobbm}. More precisely, Bona and Cascaval \cite{bona} obtained local well-posedness in Sobolev space $H^1$ and Mugnolo and Rault \cite{Mugnolobbm} showed the existence of traveling waves for the BBM equation on graphs. Using a different approach Ammari  and Crépeau  \cite{AmCr1} derived results of well-posedness and, also, stabilization for the Benjamin-Bona-Mahony equation in a star-shaped network with bounded edges. 

In this aspect, regarding control theory and inverse problems,  let us cite some previous works. Ignat \textit{et al.} in \cite{Ignat2011} worked on the inverse problem for the heat equation and the Schr\"odinger equation on a tree. Later on, Baudouin and Yamamoto \cite{Baudouin} proposed  a unified and simpler method to study the inverse problem of determining a coefficient. Results of stabilization and boundary controllability for KdV equation on star-shaped graphs was also proved by Ammari and Crépeau  \cite{AmCr} and Cerpa \textit{et al.} \cite{Cerpa,Cerpa1}. 

We caution that this is only a small sample of the extant work on graphs structure for partial differential equations.

\subsection{Presentation of the model} Let us now present the equation that we will study in this paper. The fourth-order nonlinear Schr\"odinger (4NLS) equation or biharmonic cubic nonlinear Schr\"odinger equation 
\begin{equation}
\label{fourtha}
i\partial_tu +\partial_x^2u-\partial_x^4u=\lambda |u|^2u,
\end{equation}
have been introduced by Karpman \cite{Karpman} and Karpman and Shagalov \cite{KarSha} to take into account the role of small fourth-order dispersion terms in the propagation of intense laser beams in a bulk medium with Kerr nonlinearity. Equation \eqref{fourtha} arises in many scientific fields such as quantum mechanics, nonlinear optics and plasma physics, and has been intensively studied with fruitful references (see \cite{Ben,CuiGuo,Karpman,Paus,Paus1} and references therein).

In the past twenty years such 4NLS have been deeply studied from different mathematical points of view. For example, Fibich \textit{et al.} \cite{FiIlPa} worked various properties of the equation in the subcritical regime, with part of their analysis relying on very interesting numerical developments.  The well-posedness and existence of the solutions has been shown (see, for instance, \cite{Paus,Paus1,tsutsumi,Tzvetkov}) by means of the energy method, harmonic analysis, etc.

Recently, in \cite{CaCa}, the first and the second authors worked with equation \eqref{fourtha} with the purpose to obtain controllability results. More precisely, they proved that on torus $\mathbb{T}$,  the solution of the associated linear system \eqref{fourtha} is globally exponential stable,  by using certain properties of propagation of compactness and regularity in Bourgain spaces. This property, together with the local exact controllability, ensures that fourth order nonlinear Schr\"odinger is globally exactly controllable, we suggest the reader to see \cite{CaCa} for more details.

Considering another domain instead of the torus $\mathbb{T}$, the authors, in \cite{CaCaGa},  considered the cubic fourth order  Schr\"odinger equation on the right half-line 
\begin{equation}\label{fourthc}
\begin{cases}
i\pt u +\gamma\px^4 u + \lambda |u|^2u=0, & (t,x)\in (0,T)  \times(0,\infty),\\
u(0,x)=u_0(x),                                   & x\in(0,\infty),\\
u(t,0)=f(t),\ u_x(t,0)=g(t)& t\in(0,T),
\end{cases}
\end{equation}
for $\gamma, \lambda\in\mathbb{R}$. When $\gamma\lambda< 0$ system \eqref{fourthc} is so-called focusing otherwise, that is, $\gamma\lambda > 0$, is called defocusing. In \cite{CaCaGa}, Capistrano-Filho \textit{et al.} consider $\gamma=-1$ and suitable choices of $f(t)$ and $g(t)$ in the equation \eqref{fourthc}, precisely, by assuming $$( u_0,f,g)\in H^s(\mathbb{R}^+)\times H^{\frac{2s+3}{8}}(\mathbb{R}^+)\times H^{\frac{2s+1}{8}}(\mathbb{R}^+),$$ they obtained local well-posedness on the Sobolev spaces $H^s(\R^+)$ for $s\in [0,\frac12)$. For $s > 1/2$, by the Sobolev embedding and the energy method one can easily show the local well-posedness in $H^s(\R^+)$,  giving a starting point of the study on a series of problems concerning of the Biharmonic NLS on bounded domains or star graphs.

Due these results presented in this recent work, naturally, we should see what happens for the system \eqref{fourthc}  in star graph structure given by $N$ unbounded edges $(0,\infty)$ connected with a common vertex $\nu=0$, where a function on the  graph $\mathcal{G}$ is a vector $u(t,x)=(u_1(t,x),u_2(t,x),...,u_N(t,x))$. Thus, let us consider the fourth order nonlinear Sch\"odinger equation on $\mathcal{G}$, given by 
\begin{equation}\label{graph}
\begin{cases}
i\pt u_j - \px^4 u_j + \lambda |u_j|^2u_j=0, & (t,x)\in (0,T)  \times(0,\infty),\ j=1,2, ..., N\\
u_j(0,x)=u_{j0}(x),                                   & x\in(0,\infty),\ j=1,2, ..., N\\
\end{cases}
\end{equation}
with initial conditions $(u_1(0,x),u_2(0,x), ...,u_N(0,x))$. 
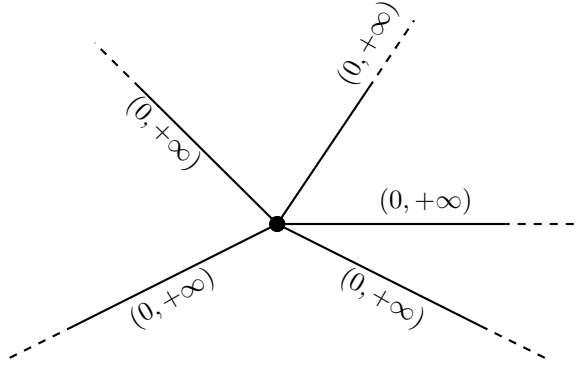
\begin{figure}[htp]
	\centering 
	\begin{tikzpicture}[scale=3]
	%\draw [thick, dashed] (-1.3,0)--(-1,0);
	%\draw[thick](-1,0)--(0,0);
	%\node at (-0.65,0.1){$(0, +\infty)$};
	
	\draw [thick, dashed] (1.3,0)--(1,0);
	\draw[thick](1,0)--(0,0);
	\node at (0.65,0.1){$(0,+\infty)$};

	\draw[thick](0,0)--(-0.6,0.6);
	\draw [thick, dashed] (-0.6,0.6)--(-0.8,0.8);
	\node at (-0.5,0.4)[rotate=-45]{$(0,+\infty)$};

	\draw[thick](0,0)--(-0.89,-0.45);
	\node at (-0.46,-0.33)[rotate=30]{$(0, +\infty)$};
	\draw [thick, dashed] (-0.89,-0.45)--(-1.19,-0.6);
	\fill (0,0)  circle[radius=1pt];

	\draw[thick](0,0)--(0.4,0.6);
	\node at (0.4,0.8)[rotate=60]{$(0,+\infty)$};
	\draw [thick, dashed] (0.4,0.6)--(0.6,0.9);
	\draw[thick](0,0)--(0.89,-0.45);
	\node at (0.46,-0.33)[rotate=-30]{$(0,+\infty)$};
	\draw [thick, dashed] (0.89,-0.45)--(1.19,-0.6);
	\fill (0,0)  circle[radius=1pt];
	\end{tikzpicture}
	\caption{Star graph with $5$ edges}
	\label{2figure}
\end{figure}

Therefore, the first following natural question arise.
%&\sum_{j=1}^N\int_0^{\infty}|u_j(x,0)|^2dx\\

\vspace{0.2cm}
\noindent\textbf{Problem $\mathcal{A}$}: Which are the boundary conditions that we can impose, at least mathematically acceptable, to ensure the well-posedness result for the system \eqref{graph}?

\subsection{Choosing the boundary conditions and main result} 

We are interested to prove the well-posedness of \eqref{graph} with appropriate boundary condition, more precisely, we will solve \eqref{graph} with the following  boundary conditions:
\begin{equation}\label{boundary1}
\text{Type $\mathcal{A}$: }
\begin{cases}
\partial_x^k u_1(t,0)=\partial_x^k u_2(t,0)=\cdots=\partial_x^k u_N(t,0),\ k=0,1 & t\in(0,T),\\
\sum_{j=1}^N\partial_x^ku_j(t,0)=0,\ k=2,3\ & t\in(0,T),
\end{cases}
\end{equation}
\begin{equation}\label{boundary2}
\text{Type $\mathcal{B}$: }\begin{cases}
\partial_x^k u_1(t,0)=\partial_x^k u_2(t,0)=\cdots=	\partial_x^k u_N(t,0),\ k=2,3 & t\in(0,T),\\
\sum_{j=1}^N\partial_x^ku_j(t,0)=0,\ k=0,1\ & t\in(0,T),
\end{cases}
\end{equation}
and
\begin{equation}\label{boundary3}
\text{Type $\mathcal{C}$: }\begin{cases}
\partial_x^k u_1(t,0)=\partial_x^k u_2(t,0)=\cdots=\partial_x^k u_N(t,0),\ k=0,3 & t\in(0,T),\\
\sum_{j=1}^N\partial_x^ku_j(t,0)=0,\ k=1,2\ & t\in(0,T).
\end{cases}
\end{equation}

These boundary conditions are motivated by the conservation of the mass. Let us denote the mass as
\begin{equation*}
\begin{split}
E(u_1(t,x)&,u_2(t,x),\cdots, u_N(t,x))=\frac12\sum_{j=1}^N\int_0^{\infty}|u_j(t,x)|^2dx.
\end{split}
\end{equation*}
Multiplying \eqref{graph} by $\overline {u}_j$, taking the imaginary part, integrating by parts and using the initial conditions of \eqref{graph}, we can obtain the most basic energy identity, namely the $L^2$--energy, satisfying 
\begin{equation}\label{energy}
\begin{split}
E(u_1(T,x),u_2(T,x),\cdots, u_N(T,x))=&- \sum_{j=1}^N \int_0^T\text{Im}(\partial_x^3u_j(t,0)\overline{u}_j(t,0)) dt\\ &+\sum_{j=1}^N\int_0^T\text{Im}(\partial_x^2u_j(t,0)\partial_x\overline{u}_j(t,0)) dt\\
&-E(u_1(0,x),u_2(0,x),\cdots, u_N(0,x)).
\end{split}
\end{equation}
Analyzing \eqref{energy}, we are interested in boundary conditions to the Cauchy problem \eqref{graph} such that the right hand side of \eqref{energy} ensures the conservation of the mass. In this sense, the boundary conditions \eqref{boundary1}, \eqref{boundary2} and \eqref{boundary3} are appropriated. Assuming one of the boundary conditions  \eqref{boundary1}, \eqref{boundary2} or \eqref{boundary3} the mass is conserved, i.e.,
\begin{equation*}
E(u_1(t,x),u_2(t,x),\cdots, u_N(t,x))=E(u_1(0,x),u_2(0,x),\cdots, u_N(0,x)).
\end{equation*}

It is important to point out that the boundary conditions of types $\mathcal A$, $\mathcal{B}$ or $\mathcal{C}$ are coherent with the study of biharmonic operator on $L^2(\mathcal{G})$. More precisely, a simple calculation proves that the  biharmonic operator $$B:=i\partial_x^4: \mathcal{D}(B_i)\subset L^2(\mathcal G)\rightarrow L^2(\mathcal{G}), \ \ i=1,2,3,$$ with the following domains  
\begin{equation*}
\begin{split}
\mathcal{D}(B_1)=\{H^4(\mathcal{G});&\ 
\partial_x^k u_1(0)=\partial_x^k u_2(0)=\cdots=\partial_x^k u_N(0),\ k=0,1\\&
\sum_{j=1}^N\partial_x^ku_j(0)=0,\ k=2,3\ \},
\end{split} 
\end{equation*}
\begin{equation*}
\begin{split}
\mathcal{D}(B_2)=\{H^4(\mathcal{G});&\ 
\partial_x^k u_1(0)=\partial_x^k u_2(0)=\cdots=\partial_x^k u_N(0),\ k=2,3\\&
\sum_{j=1}^N\partial_x^ku_j(0)=0,\ k=0,1\ \},
\end{split} 
\end{equation*}
or
\begin{equation*}
\begin{split}
\mathcal{D}(B_3)=\{H^4(\mathcal{G});&\ 
\partial_x^k u_1(0)=\partial_x^k u_2(0)=\cdots=\partial_x^k u_N(0),\ k=0,3\\&
\sum_{j=1}^N\partial_x^ku_j(0)=0,\ k=1,2\ \},
\end{split} 
\end{equation*}
is self-adjoint. Then, by Stone's Theorem (see e.g. \cite{Ca}), $B$ generates a linear group, denoted by $e^{it\partial_x^4}$ that solves the linear problem
\begin{equation*}
\begin{cases}
\partial_tu(x,t)=i\partial_x^4 u(x,t),\\
u(0,x)=u_0\in \mathcal{D}(B_i),\\
 u \in C(\R;\mathcal{D}(B_i))\cap C^1(\R;L^2(\mathcal G))\ \ i=1,2 \text{ or } 3.
\end{cases}
\end{equation*}
By using the Duhamel formula together with the fact that $H^4(\mathcal{G})$ is a Banach algebra it is possible to show that problem \eqref{graph} is well posed in high regularity, precisely, in $\mathcal{D}(B_i)$, $i=1,2$ or $3$.

\begin{remarks*} 
The following remarks are now in order.
\begin{itemize}
\item Considering the Schr\"odinger equation on a star graph $\mathcal{G}$, the vertex condition Type $\mathcal A$, when restrict on the cases $k=0$ and $k=1$, coincides with the classical Kirchhoff vertex condition. For this system, these conditions are rather natural in the context of water (and other fluids) waves, corresponding to continuity of the flow and flux balance.
\item In this direction, we cite a very interesting  work of Gregorio and Mugnolo \cite{GM} that treated the bi-laplacian on star graphs and trees with bounded edges, more precisely, they given a characterization of complete graphs in terms of the Markovian property of the semigroup generated by $\mathcal{L}^2 (\mathcal G)$, the square of the discrete Laplacian acting on a connected discrete graph $\mathcal G$. For a complete picture about  star graphs in unbounded edges, in the context of the Airy equation, we cite the work of Mugnolo \textit{et al.} \cite{MNS}.
\end{itemize}
\end{remarks*}

Therefore, this work gives an answer for the Problem $\mathcal{A}$, in a star graph structure, when the boundary conditions  \eqref{boundary1}, \eqref{boundary2} or \eqref{boundary3} are considered. This problem was left as an open problem in \cite{CaCaGa}. Before to enunciate the principal result of this work, we will denote the classical Sobolev space on the star graph $\mathcal{G}$ by
$$
H^{s}(\mathcal{G})=  \bigoplus_{i =1}^N H^{s}(0,+\infty),\quad \text{for}\,\, s\geq0.
$$
With this notation, the main result of this work can be read as follows.
\begin{theorem}\label{theorem1}
	Let $s \in [0,\frac12) $. For given initial-boundary data $(u_{10},u_{20},..., u_{N0})\in H^{s}(\mathcal{G})$ satisfying type $\mathcal{A}$, $\mathcal{B}$ or $\mathcal{C}$ vertex conditions, there exist a positive time $T$ depending on $\sum_{j=1}^N\|u_{j0}\|_{H^s(\mathbb{R}^+)}$  and a distributional solution $u=(u_i)_{i=1}^{N}(t,x) \in C((0 , T);H^s(\mathcal{G}))$ to \eqref{graph}-\eqref{boundary1} (or \eqref{graph}-\eqref{boundary2} or \eqref{graph}-\eqref{boundary3}) satisfying
\[\begin{aligned}
u_j \in C&\bigl(\mathbb{R}^{+};\; H^{\frac{2s+3}{8}}(0,T)\bigr) \cap X^{s,b}(\R^+\times (0,T)) , \\
&\partial_xu_j\in C\bigl(\R^+;\; H^{\frac{2s+1}{8}}(0,T)\bigr), \\
& \partial_x^2u_j\in C\bigl(\R^+;\; H^{\frac{2s-1}{8}}(0,T)\bigr)
\end{aligned}\]
and
\[\begin{aligned}
& \partial_x^3u_j\in C\bigl(\R^+;\; H^{\frac{2s-3}{8}}(0,T)\bigr),
\end{aligned}\]
for some $b(s) < \frac12$ and $j=1,2,..., N$. Moreover, the map $(u_{10},u_{20},...,u_{N0})\longmapsto u$ is locally Lipschitz continuous from $H^s(\mathcal{G})$ to  $C\big((0,T);\,H^s(\mathcal{G})\big)$.
\end{theorem}

\subsection{Heuristic of the paper and further comments} In this work we prove the existence of solution to the problem \eqref{graph}-\eqref{boundary1} (or \eqref{graph}-\eqref{boundary2} or \eqref{graph}-\eqref{boundary3}) on  star graph structure $\mathcal{G}$ composed by $N$ unbounded edges. The proof of Theorem \ref{theorem1} will be divided in several steps.  Initially, we recast the partial differential equation in each edge for a full line with a forcing term, more precisely
\begin{equation}\label{extended}
\begin{cases}
i\pt u_j - \px^4 u_j + \lambda |u_j|^2u_j=\mathcal{T}_{1j}(x) h_{1j}(t)+\mathcal{T}_{2j}(x) h_{2j}(t), & (t,x)\in (0,T)  \times \R,\ j=1,2, ..., N\\
u_j(0,x)=\widetilde u_{j0}(x),                                   & x\in\R,\ j=1,2, ..., N\\
\end{cases}
\end{equation}
where $\mathcal{T}_{1j}$ and $\mathcal{T}_{2j}$ $(j=1,2,...,N)$ are distributions supported in the negative half-line $(-\infty,0)$; the
boundary forcing functions $h_{1j},h_{2j}$  $(j=1,2,...,N)$ are selected to ensure that the vertex conditions are satisfied and $\widetilde u_{j0}(x)$ are extensions of $u_{j0}$ $(j=1,2,...,N)$ on full line satisfying
$$\|\widetilde u_{j0}\|_{H^s(\R)}\leq 2 \| u_{j0}\|_{H^s(\R^+)}.$$
%The boundary forcing functions $h_{1j}, h_{2j}$, for $j=1,2,...,N$, are selected to ensure that the vertex conditions are satisfied, by inverting a Riemann-Liouville fractional integral. Here, we will consider the distributions $\mathcal{T}_{kj}=\frac{x_{-}^{-1+\lambda_{k,j}}}{\Gamma(\lambda_{k,j})}$ defined by
%\begin{equation}
%\left\langle \frac{x_-^{\lambda-1}}{\Gamma(\lambda)},\phi\right\rangle=\int_0^{+\infty}\frac{(-x)^{\lambda-1}}{\Gamma(\lambda)}\phi(x)dx,
%\end{equation}
%for $\text{Re}\ \lambda>0$ and for $\text{Re}\ \lambda<0$ we use the formula $\frac{x_-^{\lambda-1}}{\Gamma(\lambda)}=\frac{d^k}{dx}\frac{x_-^{\lambda+k-1}}{\Gamma(\lambda+k)}$ to define the distribution $\frac{x_-^{\lambda-1}}{\Gamma(\lambda)}$.  The crucial point here is to choices, appropriately,  the parameters $\lambda_{k,j}$ and the functions $h_{k,j}$, for $k=1,2$ and $j=1,2,...,N$, that will depend on the regularity index $s$ and the boundary condition considered. This will be chosen by putting all boundary conditions in term of matrices.
Upon constructing the solution $\tilde u=(\tilde u_1,\tilde u_2,...,\tilde u_N)$ of \eqref{extended}, we obtain
the solution $u=(u_1,u_2,...,u_N)$ of problem \eqref{graph} with appropriate boundary condition, by restriction, as 
$$u=u(x,t)=(u_1,u_2,...,u_N)\big|_{x\in \mathcal{G}, t\in (0,T)}:=(u_1|_{x\in \R^+, t\in (0,T)},u_2|_{x\in \R^+, t\in (0,T)},...,u_N|_{x\in \R^+, t\in (0,T)}).$$
Secondly,  the solution of forced Cauchy problem \eqref{extended} satisfying the vertex types $\mathcal{A}$, $\mathcal{B}$ or $\mathcal{C}$,  is constructed using the classical Fourier restriction method due Bourgain \cite{Bourgain1993}. Finally, a fixed point argument ensures the proof of the Theorem \ref{theorem1}.

\vspace{0.2cm}

We present some comments about the relevance of the method used in this manuscript.
\begin{itemize}
\item[i.] It is important to point out that,  in our knowledge, this work is the first one in a star graphs structure $\mathcal{G}$ composed by $N$ unbounded edges by using boundary forcing operator approach introduced first by Colliander and Kenig \cite{CK} and improved by Holmer \cite{Holmerkdv}.

\vspace{0.1cm}

\item[ii.] The graph structure of this article is more complex than proposed in \cite{Cav} in the following sense:
To treat the extended vectorial integral equation that solves system \eqref{graph}, considering $ N $ unbounded edges with appropriated vertex conditions, is more delicate since the matrices associated with this problem have $2N$--order (see Section \ref{sec4}).

%\vspace{0.1cm}

\item[iii.] A more delicate question concerns here is the local well-posedness for the Cauchy problem \eqref{graph} in low regularity. To do this we need to use a dispersive approach instead of Semigroup theory, where the principal difficulty is to use the restriction Fourier method in the context of star graphs. This motivates us to solve the problem \eqref{graph} by using this approach, since the Semigroup theory does not guarantee the lower regularity to solutions of \eqref{graph}.
%\item[iii.] We decided do not explore the local well-posedness result on higher regularity $H^s(\mathcal{G})$, i.e., for $s>\frac12$, because this analysis is more direct since Sobolev spaces is a Banach algebra for $s>\frac12$ and it is not required the use of Bourgain spaces. Additionally, instead to use boundary forcing operator, for $s>\frac12$, we can use Semigroup theory to prove the result in higher regularity. 
\end{itemize}

\subsection{Organization of the article} 
To end our introduction, we present the outline of the manuscript. Section \ref{sec2} is devoted to present the notations, more precisely, the Sobolev spaces, the Bourgain spaces, the Riemann-Liouville fractional integral operator and the Duhamel boundary forcing operator associated of (4NLS), which are paramount to prove the main result of the article. In the section \ref{sec3}, we  will give an overview of the main estimates proved by the authors in \cite{CaCaGa}.  With these two sections in hand, we are able to prove Theorem \ref{theorem1}, in several steps, in the Section \ref{sec4}. The Section \ref{sec5} is devoted to prove an auxiliary lemma, which one is used in the proof of the main result of the article, namely, Theorem \ref{theorem1}. Finally, at the end of the work, we present an Appendix \ref{apenB}, which will we given a sketch of the proof of Theorem \ref{theorem1} with vertex conditions types $\mathcal{B}$ and $\mathcal{C}$. 

\section{Preliminaries}\label{sec2}
This section is devoted to presenting the main notations, introducing the functions spaces used in this work and the Duhamel boundary forcing operator associated with the fourth order linear Schr\"odinger equation.
\subsection{Notations}
Let us fix a cut-off function $\psi(t):=\psi$ such that $\psi \in C_0^{\infty}(\mathbb{R})$, $0 \le \psi \le1$ and defined by
\begin{equation}\label{eq:cutoff}
\psi \equiv 1 \; \mbox{ on } \; [0,1], \quad \psi \equiv 0, \;\text{for } |t| \ge 2, 
\end{equation}
and, for $T>0$, we denote $\psi_{T}(t)=\frac{1}{T}\psi(\frac{t}{T})$.

Now, for $s\geq 0$, define the homogeneous $L^2$-based Sobolev spaces $\dot{H}^s=\dot{H}^s(\R)$ by natural norm $\|\phi\|_{\dot{H}^s}=\||\xi|^s\hat{\psi}(\xi)\|_{L^2_{\xi}}$ and the $L^2$-based inhomogeneous Sobolev spaces $H^s=H^s(\mathbb{R})$ by the norm $\|\phi\|_{\dot{H}^s}=\|(1+|\xi|^2)^{\frac{s}{2}}\hat{\psi}(\xi)\|_{L^2_{\xi}}$, where $\hat{\psi}$ denotes the Fourier transform of $\psi$.  The function $f$ belongs to $H^s(\mathbb{R}^+)$, if there exists $F \in H^s(\R)$ such that $f(x)=F(x)$ for $x>0$, in this case we set $$\|f\|_{H^s(\mathbb{R}^+)}=\inf_{F}\|F\|_{H^{s}(\mathbb{R})}.$$ On the other hand, for $s \in \R$, $f \in H_0^s(\mathbb{R}^+)$ provided that there exists $F \in H^s(\R)$ such that $F$ is the extension of $f$ on $\R$ and $F(x) = 0$ for $x<0$. In this case, we set $\norm{f}_{H^s_0(\R^+)} = \inf_{F} \norm{F}_{H^s(\R)}$. 
For $s<0$, we define $H^s(\mathbb{R}^+)$ as the dual space of $H_0^{-s}(\mathbb{R}^+)$. It is well known that $H_0^s(\R^+)=H^{s}(\R^+)$ for $-\frac12<s<\frac12$.

Finally, the sets $C_0^{\infty}(\mathbb{R}^+)=\{f\in C^{\infty}(\mathbb{R});\, \supp f \subset [0,\infty)\}$ and $C_{0,c}^{\infty}(\mathbb{R}^+)$ are defined as the subset of $C_0^{\infty}(\mathbb{R}^+)$, whose members have a compact support on $(0,\infty)$. We remark that $C_{0,c}^{\infty}(\mathbb{R}^+)$ is dense in $H_0^s(\mathbb{R}^+)$ for all $s\in \mathbb{R}$.

\subsection{Solution spaces}\label{sec:sol space}
Consider $f \in \Sch (\R^2)$, let us denote by $\wt{f}$ or $\ft (f)$ the Fourier transform of $f$ with respect to both spatial and time variables\[\wt{f}(\tau , \xi)=\int _{\R ^2} e^{-ix\xi}e^{-it\tau}f(t,x) \;dxdt .\]Moreover, we use $\ft_x$ and $\ft_t$ to denote the Fourier transform with respect to space and time variables, respectively (also use $\,\, \wh{\;}\,\,$ for both cases). 
%\marginnote{\color{blue}Looking the calculations!!! I guess that we have to put $\frac{1}{2\pi}$ in the fourier transform definition}

In the 90's Bourgain \cite{Bourgain1993} established a form to show the well-posedness of some classes of dispersive systems. Precisely, on the Sobolev spaces $H^s$, for smaller values of $s$, with these new spaces, Bourgain showed a smoothing property more suitable for solutions of these classes of dispersive equations.

In our case, considering $s,b\in \mathbb{R}$ we present below the Bourgain spaces $X^{s,b}$ associated to the linear system of  \eqref{graph}. The space will be a completion of $\Sch'(\mathbb{R}^2)$ under the norm\[\norm{f}_{X^{s,b}}^2 = \int_{\R^2} \bra{\xi}^{2s}\bra{\tau +\xi^4}^{2b}|\wt{f}(\tau,\xi)|^2 \; d\xi d\tau, \]where $\bra{\cdot} = (1+|\cdot|^2)^{1/2}$.

It is important to note that $X^{s,b}$--space, with $b> \frac12$, is well-adapted to study the IVP of dispersive equations. However, in the study of IBVP, the standard argument cannot be applied directly. This is due to the lack of hidden  regularity, more precisely, the control of (derivatives) time trace norms of the Duhamel parts requires to work in $X^{s,b}-$type spaces for $b<\frac12$, since the full regularity range cannot be covered (see Lemma \ref{duhamel} inequality \eqref{derivative1}). 

Considering the space denoted by $Z^{s,b}$ with the following norm
\[\norm{f}_{Z^{s,b}(\R^2) }= \sup_{t \in \R} \norm{f(t,\cdot)}_{H^s(\R)} + \sum_{j=0}^{3}\sup_{x \in \R} \norm{\px^jf(\cdot,x)}_{H^{\frac{2s+3-2j}{8}}(\R)} + \norm{f}_{X^{s,b}} ,\]
our goal is to find solutions of the Cauchy problem \eqref{graph}. 

Here, we will consider the spatial and time restricted space of $Z^{s,b}(\R^2)$ defined in the standard way as follows
\[Z^{s,b}((0 , T)\times \R^+) = Z^{s,b} \Big|_{(0 , T)\times \R^+}\]
equipped with the norm
\[\norm{f}_{Z^{s,b}((0 , T)\times \R^+) } = \inf\limits_{g \in Z^{s,b}} \set{\norm{g}_{Z^{s,b}} : g(t,x) = f(t,x) \; \mbox{ on } \; (0,T) \times \R^+}. \]

\subsection{Riemann-Liouville fractional integral} Before beginning our study of the Cauchy problem \eqref{graph}, in this subsection, we just give a brief summary of the Riemann-Liouville fractional integral operator to make the work complete. We suggest \cite{CaCaGa,CK,Holmerkdv} for the reader to see the proofs and more details.

Consider the function $t_+$ defined in the following way
\[t_+ = t \quad \mbox{if} \quad t > 0, \qquad t_+ = 0  \quad \mbox{if} \quad t \le 0.\]
The tempered distribution $\frac{t_+^{\alpha-1}}{\Gamma(\alpha)}$ is defined as a locally integrable function by
\begin{equation*}
	\left \langle \frac{t_+^{\alpha-1}}{\Gamma(\alpha)},\ f \right \rangle=\frac{1}{\Gamma(\alpha)}\int_0^{\infty} t^{\alpha-1}f(t)dt,
\end{equation*}
 for Re $\alpha>0$. By integrating by parts, we have that 
\begin{equation}\label{gamma}
	\frac{t_+^{\alpha-1}}{\Gamma(\alpha)}=\partial_t^k\left( \frac{t_+^{\alpha+k-1}}{\Gamma(\alpha+k)}\right),
\end{equation}
for all $k\in\mathbb{N}$. Expression  \eqref{gamma} allows to extend the definition of $\frac{t_+^{\alpha-1}}{\Gamma(\alpha)}$, in the sense of distributions, to all $\alpha \in \mathbb{C}$. 
For $f\in C_0^{\infty}(\mathbb{R}^+)$, define $	\mathcal{I}_{\alpha}f$ as
\begin{equation*}
	\mathcal{I}_{\alpha}f=\frac{t_+^{\alpha-1}}{\Gamma(\alpha)}*f.
\end{equation*}
Thus, for $\mbox{Re }\alpha>0$, follows that
\begin{equation}\label{eq:IO}
	\mathcal{I}_{\alpha}f(t)=\frac{1}{\Gamma(\alpha)}\int_0^t(t-s)^{\alpha-1}f(s) \; ds
\end{equation}
and the following properties easily holds $$\mathcal{I}_0f=f,\quad \mathcal{I}_1f(t)=\int_0^tf(s)ds, \quad\mathcal{I}_{-1}f=f '\quad \text{and} \quad \mathcal{I}_{\alpha}\mathcal{I}_{\beta}=\mathcal{I}_{\alpha+\beta}.$$ 

\subsection{Duhamel boundary forcing operator}\label{sec:Duhamel boundary forcing operator}
Now, we present the Duhamel boundary forcing operator, which was introduced by Colliander and Kenig \cite{CK}, in order to construct the solution to
\begin{equation}\label{eq:linear 4kdv}
i\partial_t u-\partial_x^4u=0.
\end{equation}
For details about this subsection and for a well exposition about this topic, the authors suggest the following references \cite{CaCaGa,Cavalcante,CC,HolmerNLS} .

Following \cite{CaCaGa}, let us consider the oscillatory integral by
\begin{equation}\label{eq:oscil}
B(x)=\frac{1}{2\pi}\int_{\R}e^{ix\xi}e^{-i \xi^4} \; d\xi,
\end{equation} 
which one is the key to define the Duhamel boundary forcing operator. A change of variable and contour proves that $B(0) =-\frac{i^{7/4}}{\pi}\Gamma\left( \frac54 \right).$ We will denote
\begin{equation}\label{eq:M}
M = \frac{1}{ B(0)\Gamma(3/4)}.
\end{equation} 
For $f\in C_0^{\infty}(\mathbb{R}^+)$, define the boundary forcing operator $\mathcal{L}^0$ (of order $0$) as
\begin{equation}\label{eq:BFO}	
\mathcal{L}^0f(t,x):=M\int_0^te^{i(t-t')\partial_x^4}\delta_0(x)\mathcal{I}_{-\frac34}f(t')dt',
\end{equation}
where $e^{it\partial_x^4}$ denotes the group associated to \eqref{eq:linear 4kdv} given by 
\begin{align*}
e^{it\partial_x^4}\psi(x)=\frac{1}{2\pi}\int_{\mathbb{R}} e^{i x \xi}e^{-it\xi^4}\hat{\psi}(\xi)d\xi.
\end{align*}
By using the following properties of  the convolution operator (for $k=1$)
\begin{equation}\label{conv}
\partial_x^k (f * g) = (\partial_x^lf) * g = f * (\partial^k_xg), \quad k\in\mathbb{N},
\end{equation} 
and the integration by parts in $t'$ of \eqref{eq:BFO}, we get that
\begin{equation}\label{eq:BFO1}
i\mathcal{L}^0(\partial_t f)(t,x) = iM\delta_0(x)\mathcal{I}_{-\frac34}f(t) + \partial_x^4\mathcal{L}^0f(t,x).
\end{equation}
Using \eqref{eq:oscil} and, again, by change of variable, we have
\begin{equation}\label{forcing}
\begin{split}
	\mathcal{L}^0f(t,x)&=M\int_0^te^{i(t-t')\partial_x^4}\delta_0(x)\mathcal{I}_{-\frac34}f(t')dt'\\
	&= M \int_0^t B\left(\frac{x}{(t-t')^{1/4}}\right)\frac{\mathcal{I}_{-\frac34}f(t')}{(t-t')^{1/4}}dt'.
\end{split}	
\end{equation}

We are now generalize the boundary forcing operator \eqref{eq:BFO}. For $\mbox{Re}\ \lambda > -4$ and given $g\in C_0^{\infty}(\mathbb{R}^+)$, we define
\begin{equation}\label{eq:BFOC}
\mathcal{L}^{\lambda}g(t,x)=\left[\frac{x_{-}^{\lambda-1}}{\Gamma(\lambda)}*\mathcal{L}^0\big(\mathcal{I}_{-\frac{\lambda}{4}}g\big)(t, \cdot)   \right](x),
\end{equation}
where $*$ denotes the convolution operator and
$\frac{x_{-}^{\lambda-1}}{\Gamma(\lambda)}=\frac{(-x)_+^{\lambda-1}}{\Gamma(\lambda)}$.
In particular, for Re $\lambda>0$, we have
 \begin{equation}\label{forcing2}
\mathcal{L}^{\lambda}g(t,x)=\frac{1}{\Gamma(\lambda)}\int_{x}^{\infty}(y-x)^{\lambda-1}\mathcal{L}^0\big(\mathcal{I}_{-\frac{\lambda}{4}}g\big)(t,y)dy.
\end{equation}
By using the property \eqref{conv}, for $k=4$, and \eqref{eq:BFO1} give us
\begin{equation}\label{classe22}
\begin{aligned}
\mathcal{L}^{\lambda}g(t,x)&=\left[\frac{x_{-}^{(\lambda+4)-1}}{\Gamma(\lambda+4)}*\px^4\mathcal{L}^0\big(\mathcal{I}_{-\frac{\lambda}{4}}g\big)(t, \cdot)   \right](x)\\
&=iM\frac{x_{-}^{(\lambda+4)-1}}{\Gamma(\lambda+4)}\mathcal{I}_{-\frac{3}{4}-\frac{\lambda}{4}}g(t)+i\int_{x}^{\infty}\frac{(y-x)^{(\lambda+4)-1}}{\Gamma(\lambda+4)}\mathcal{L}^0\big(\partial_t\mathcal{I}_{-\frac{\lambda}{4}}g\big)(t,y)dy,
\end{aligned}
\end{equation}
for $\mbox{Re}\  \lambda > -4$, where $M$ is defined as \eqref{eq:M}.  From \eqref{eq:BFO1} and \eqref{eq:BFOC}, we have
\begin{equation*}
(i\partial_t-\partial_x^4)\mathcal{L}^{\lambda}g(t,x)=iM\frac{x_{-}^{\lambda-1}}{\Gamma(\lambda)}\mathcal{I}_{-\frac{3}{4}-\frac{\lambda}{4}}g(t),
\end{equation*}
in the distributional sense. 
\section{Overview of the main estimates}\label{sec3}
With all the notations and spaces defined in Section \ref{sec2}, we present now the main estimates of this work, which are paramount to prove the main result of the article. 

\subsection{Estimates for the function spaces} 
Concerning of the $X^{s,b}$ space, we have two properties which are presented in the lemma below for the functions $\psi(t)$ and $\psi_T$, defined in \eqref{eq:cutoff}. The first item can be found in  \cite[Lemma 2.11]{Tao2006} and the second one in Ginibre \textit{et al.} \cite{GTV}. The lemma can be read as follows:
\begin{lemma}\label{gvt}
Let $\psi(t)$ be a Schwartz function in time. Then, we have 
\[\norm{\psi(t)f}_{X^{s,b}} \lesssim_{\psi,b} \norm{f}_{X^{s,b}}.\]
Moreover, if $-\frac{1}{2}< b'< b\leq 0$\, or\, $0\leq b'<b<\frac{1}{2}$, $w\in X^{s,b}$ and $s\in \mathbb{R}$, thus 
	\begin{equation*}
	\|\psi_{T}w\|_{ X^{s,b'}}\leq c T^{b-b'}\|w\|_{ X^{s,b}}.
	\end{equation*}
\end{lemma}

Another result that state important properties of the Riemann-Liouville fractional integral operator is given below. The proof of this can be found in \cite[Lemmas 2.1, 5.3 and 5.4]{Holmerkdv}.
\begin{lemma}\label{lio} If $f\in C_0^{\infty}(\mathbb{R}^+)$, then $\mathcal{I}_{\alpha}f\in C_0^{\infty}(\mathbb{R}^+)$, for all $\alpha \in \mathbb{C}$. Moreover, we have the following:
\begin{itemize} 
\item[(a)] If $0\leq \mathrm{Re} \ \alpha <\infty$ and $s\in \mathbb{R}$, then $\|\mathcal{I}_{-\alpha}h\|_{H_0^s(\mathbb{R}^+)}\leq c \|h\|_{H_0^{s+\alpha}(\mathbb{R}^+)}$, where $c=c(\alpha)$.
\item[(b)]If $0\leq \mathrm{Re}\ \alpha <\infty$, $s\in \mathbb{R}$ and $\mu\in C_0^{\infty}(\mathbb{R})$, then
	$\|\mu\mathcal{I}_{\alpha}h\|_{H_0^s(\mathbb{R}^+)}\leq c \|h\|_{H_0^{s-\alpha}(\mathbb{R}^+)},$ where $c=c(\mu, \alpha)$.

\end{itemize}
\end{lemma}

\subsection{Estimates for the boundary forcing operator} Now, let us be precisely when the boundary forcing operator is continuous or discontinuous. Initially, we present the well-know properties of the spatial continuity, the decay of the $\mathcal{L}^{\lambda}g(t,x)$ and the explicit values for $\mathcal{L}^{\lambda}f(t,0)$, respectively, these results with their respective proofs can be seen in \cite[Lemmas 3.2 and 3.3]{CaCaGa}.
\begin{lemma}\label{continuous}
	Let $g\in C_0^{\infty}(\mathbb{R}^+)$ and $M$ be as in \eqref{eq:M}. Then, we have
	\begin{equation}\label{eq:relation}
\mathcal{L}^{-k}g=\partial_x^k\mathcal{L}^{0}\mathcal{I}_{\frac{k}{4}}g, \qquad k=0,1,2,3.
	\end{equation}
	Moreover, for fixed $ 0 \le t \le 1$, $\partial_x^k \mathcal{L}^0f(t,x)$, $k=0,1,2$, is continuous in $x \in \mathbb{R}$ and $\mathcal{L}^{-3}g(t,x)$ is continuous in $x \in \R \setminus \set{0}$ and has a step discontinuity  at $x=0$.
%	 For real $\lambda$ satisfying $\lambda>-3$, $\mathcal{L}^{\lambda}g(t,x)$ is continuous in $x \in\mathbb{R}$. For $-3\leq\lambda\leq 1$ and  $0\leq t\leq 1$, $\mathcal{L}^{\lambda}g(t,x)$ satisfies the following decay bounds:
%%	&|\mathcal{L}_{-}^{\lambda}g(t,x)|\leq c_{m,\lambda,g}\langle x\rangle^{-m},\; \text{for all}\quad x\leq 0 \quad \text{and} \quad m\geq0,\\ 
%	\begin{align*}
%	&|\mathcal{L}^{\lambda}g(t,x)|\leq c_{\lambda,g}\langle x\rangle^{\lambda-1}, \quad \text{for all} \quad  x\geq 0\\
%	\intertext{and}
%	&|\mathcal{L}^{\lambda}g(t,x)|\leq c_{m,\lambda,g}\langle x\rangle^{-m}, \quad \text{for all} \quad x\geq 0 \quad \text{and} \quad m\geq0.
%	\end{align*}
%%	\intertext{and} &|\mathcal{L}_{+}^{\lambda}g(t,x)|\leq c_{\lambda,g}\langle x\rangle^{\lambda-1},\quad \text{for all}\quad  x\leq 0.
\end{lemma}

\begin{lemma}\label{trace}
	 For $\mbox{Re}\ \lambda>  -4$ and $f\in C_0^{\infty}(\R^+)$, we have the following value of $\mathcal{L}^{\lambda}f(t,0)$:
	\begin{equation}\label{lr0}
	\mathcal{L}^{\lambda}f(t,0)=\frac{M}{8\ }f(t)\left(\frac{e^{-i\frac{\pi}{8}(1+3\lambda)}+e^{-i\frac{\pi}{8}(1-5\lambda)}}{\sin(\frac{1-\lambda}{4}\pi)}\right).
	\end{equation}
	\end{lemma}
\subsection{Energy and trilinear estimates}\label{sec:energy}
In the last part, we present four lemmas related to energy and  trilinear estimates for the solutions of the 4NLS equation in the Bourgain spaces $X^{s,b}$.  These results and their proofs can also be found in \cite[Section 4]{CaCaGa}.
\begin{lemma}\label{grupo}
	Let $s\in\mathbb{R}$ and $b \in \R$. If $\phi\in H^s(\mathbb{R})$, then the following estimates holds
	\begin{equation}\label{space}
	\|\psi(t)e^{it\partial_x^4}\phi(x)\|_{C_t\big(\mathbb{R};\,H_x^s(\mathbb{R})\big)}\lesssim_{\psi} \|\phi\|_{H^s(\mathbb{R})},
	\end{equation}
	\begin{equation}\label{derivative}
	\|\psi(t) \partial_x^{j}e^{it\partial_x^4}\phi(x)\|_{C_x\left(\mathbb{R};H_t^{\frac{2s+3-2j}{8}}(\mathbb{R})\right)}\lesssim_{\psi, s, j} \|\phi\|_{H^s(\mathbb{R})} \quad j\in\N
	\end{equation}
	and
	\begin{equation}\label{bourgain}
	\|\psi(t)e^{it\partial_x^4}\phi(x)\|_{X^{s,b}}\lesssim_{\psi, b}  \|\phi\|_{H^s(\mathbb{R})}.
	\end{equation}
%	\begin{itemize}
%		\item[(a)] \emph{Space traces} $$\|\psi(t)e^{t\partial_x^4}\phi(x)\|_{C_t\big(\mathbb{R};\,H_x^s(\mathbb{R})\big)}\lesssim \|\phi\|_{H^s(\mathbb{R})};$$
%		\item[(b)] \emph{Derivative time traces} 
%	$$
%		\|\psi(t) \partial_x^{j}e^{it\partial_x^4}\phi(x)\|_{C_x(\mathbb{R};H_t^{\frac{2s+3-j}{8}}(\mathbb{R}))}\lesssim_{\psi, s, j} \|\phi\|_{H^s(\mathbb{R})}, \quad j\in\{0,1,2\};
%$$
%		\item [(c)] \emph{Bourgain spaces estimates} $$\|\psi(t)e^{it\partial_x^4}\phi(x)\|_{X^{s,b}\cap D^{\alpha}}\lesssim_{\psi, b, \alpha}  \|\phi\|_{H^s(\mathbb{R})}.$$
%	\end{itemize}

%Estimates \eqref{space}, \eqref{derivative} and \eqref{bourgain} are the so called space traces, derivative time traces and Bourgain spaces estimates, respectively.
\end{lemma}

\begin{lemma}\label{duhamel}
Let $0 < b < \frac12 $ and $j=0,1,2,3$, we have the following inequalities
\begin{equation}\label{space1}
\|\psi(t)\mathcal{D}w(t,x)\|_{C\big(\mathbb{R}_t;\,H^s(\mathbb{R}_x)\big)}\lesssim \|w\|_{X^{s,-b}},
\end{equation}
for $s\in \mathbb{R}$; 
\begin{equation}\label{derivative1}
		\|\psi(t) \partial_x^j\mathcal{D}w(t,x)\|_{C\left(\mathbb{R}_x;H^{\frac{2s+3-2j}{8}}(\mathbb{R}_t)\right)}\lesssim	\|w\|_{X^{s,-b}},
\end{equation}
for $-\frac32+j<s<\frac12+j$;
\begin{equation}\label{bourgain1}
\|\psi(t) \partial_x^j\mathcal{D}w(t,x)\|_{X^{s,b}}\lesssim	\|w\|_{X^{s,-b}},
\end{equation}
for $s\in \mathbb{R}$. 

%	\begin{itemize}
%		\item[(a)] \emph{Space traces} 
%$$
%		\|\psi(t)\mathcal{D}w(t,x)\|_{C\big(\mathbb{R}_t;\,H^s(\mathbb{R}_x)\big)}\lesssim \|w\|_{X^{s,-b}};
%$$
%		\item[(b)] \emph{Time traces}
%$$
%		\|\psi(t) \partial_x^j\mathcal{D}w(t,x)\|_{C(\mathbb{R}_x;H^{\frac{2s+3-j}{8}}(\mathbb{R}_t))}\lesssim	\|w\|_{X^{s,-b}}+\|w\|_{Y^{s,-b}};
%$$
%		\item[(c)] \emph{Bourgain spaces estimates} 
%$$
%		\|\psi(t)\mathcal{D}w(t,x)\|_{X^{s,b} }\lesssim  \|w\|_{X^{s,-b}}.
%$$
%	\end{itemize}
%Estimates \eqref{space1}, \eqref{derivative1} and \eqref{bourgain1} are so called space traces, derivative time traces and Bourgain spaces estimates, respectively.
\end{lemma}

%Remark in (b) that $\|\psi(t) \partial_x^j\mathcal{D}w(t,x)\|_{C(\mathbb{R}_x;H_0^{\frac{2s+3-j}{8}}(\mathbb{R}_t^+))}$ has same bound for $s < \frac{11}{2} + j$. 
%\marginnote{\color{red} Check these calculation with the sign minus in mind!!!!!}
\begin{lemma}\label{operator}
	Let $s\in\mathbb{R}$. Then, 
\begin{itemize}
\item[(a)] For $\frac{2s-7}{2} <\lambda< \frac{1+2s}{2}$ and $\lambda< \frac12$ the following inequality holds
		\begin{equation}\label{space2}
		\|\psi(t)\mathcal{L}^{\lambda}f(t,x)\|_{C\big(\mathbb{R}_t;\,H^s(\mathbb{R}_x)\big)}\leq c \|f\|_{H_0^\frac{2s+3}{8}(\mathbb{R}^+)};
		\end{equation}
\item[(b)] For $-4+j<\lambda<1+j$, $j=0,1,2,3$, we have 
		\begin{equation}\label{eq:(b)0}
		\|\psi(t)\partial_x^j\mathcal{L}^{\lambda}f(t,x)\|_{C\big(\mathbb{R}_x;\,H_0^{\frac{2s+3-2j}{8}}(\mathbb{R}_t^+)\big)}\leq c \|f\|_{H_0^\frac{2s+3}{8}(\mathbb{R}^+)};
		\end{equation}
\item[(c)] If $ s <4-4b $, $b < \frac12$, $-5<\lambda<\frac12$ and $s+4b-2<\lambda<s+\frac12$ yields that
\begin{equation}\label{bourgain2}
\|\psi(t)\mathcal{L}^{\lambda}f(t,x)\|_{X^{s,b}}\leq c \|f\|_{H_0^\frac{2s+3}{8}(\mathbb{R}^+)}.
\end{equation}
\end{itemize}
\end{lemma}
\begin{remarks*}
Let us present two remarks.
\begin{itemize}
\item[i.] The  previous estimates are the so-called space traces, derivative time traces and Bourgain spaces estimates, respectively.
\item[ii.] We observe that in \cite{CaCaGa} was obtained \eqref{derivative}, \eqref{derivative1} and \eqref{eq:(b)0}, for $j=0$ and $j=1$, but the result for all $j$ can be obtained directly by using the fact that $$\partial_x^j\mathcal{L}^{\lambda}=\mathcal{L}^{\lambda-j}(\mathcal{I}_{-\frac j4}).$$
\end{itemize}
\end{remarks*}
To close this section, let us enunciate the trilinear estimates associated to fourth order nonlinear Schr\"odinger equation. The proof of this estimate can be found in \cite{Tzvetkov}.
\begin{lemma}\label{trilinear}
	For $s\geq 0 $, there exists $b = b(s) < 1/2$ such that 
	\begin{equation}\label{eq:bilinear1}
	\norm{u_1u_2\overline{u}_3}_{X^{s,-b}} \le c\norm{u_1}_{X^{s,b}}\norm{u_2}_{X^{s,b}}\norm{u_3}_{X^{s,b}}.
	\end{equation}
\end{lemma}

\section{Proof of the main result}\label{sec4}
The aim of this section is to prove the main result announced  in the introduction of this work, Theorem \ref{theorem1}. Here, we only consider the vertex condition \eqref{boundary1} (type $\mathcal{A}$) and  to make the proof easy to understand, we will split it in several steps which will be divided into subsections. Additionally, the discussion of vertex conditions types $\mathcal{B}$ and $\mathcal{C}$  will be presented on Appendix \ref{apenB}, at the end of the article.
\subsection{Obtaining a integral equation in $\bigoplus_{i =1}^N\R$}\label{proof}
In this first step, we are interested in finding an extended vectorial integral equation posed in  $\bigoplus_{i =1}^N\R$, such that the restrictions of this equation on $\mathcal{G}$  will solve, in the sense of distributions, the  following Cauchy problem
\begin{equation}\label{grapha}
\begin{cases}
i\pt u_j - \px^4 u_j + \lambda |u_j|^2u_j=0, & (t,x)\in (0,T)  \times(0,\infty),\ j=1,2, ..., N\\
u_j(0,x)=u_{j0}(x),                                   & x\in(0,\infty),
\end{cases}
\end{equation}
with initial conditions $(u_{10},u_{20}, ...,u_{N0})\in H^{s}(\mathcal{G})$. Let us begin rewriting the Type $\mathcal{A}$ vertex conditions \eqref{boundary1}
in terms of matrices. In this way, note that \eqref{boundary1} is equivalent to
 \begin{equation*}
\begin{cases}
\partial_x^k u_1(t,0)=\partial_x^k u_2(t,0),\,\, \partial_x^k u_2(t,0)=\partial_x^k u_3(t,0), \cdots, \partial_x^k u_{N-1}(t,0)=\partial_x^k u_{N}(t,0), &  k=0,1, \quad t\in(0,T)\\
\sum_{j=1}^N\partial_x^ku_j(t,0)=0, & k=2,3, \quad t\in(0,T).
\end{cases}
\end{equation*}
Thus, we consider the following matrices

\begin{equation}\label{m1}
\left[ C_1 \right ]_{2N \times N}\left[\begin{array}{c}
u_1(t,0)\\
u_2(t,0) \\
u_3(t,0) \\
\vdots\\
u_{N-1}(t,0) \\
u_N(t,0)  
\end{array}\right]_{N\times 1}=\left[\begin{array}{c}
0\\
0\\
\vdots\\
0  
\end{array}\right]_{2N\times 1} ;\quad 
\left[ C_2 \right ]_{2N \times N}\left[\begin{array}{c}
\partial_x u_1(t,0)\\
\partial_xu_2(t,0) \\
\partial_xu_3(t,0) \\
\vdots\\
\partial_xu_{N-1}(t,0) \\
\partial_xu_N(t,0)  
\end{array}\right]_{N\times 1}=\left[\begin{array}{c}
0\\
0\\
\vdots\\
0  
\end{array}\right]_{2N\times 1}
\end{equation}
and
\begin{equation}\label{m3}
\left[ C_3 \right ]_{2N \times N}\left[\begin{array}{c}
\partial_x^2 u_1(t,0)\\
\partial_x^2u_2(t,0) \\
\partial_x^2u_3(t,0) \\
\vdots\\
\partial_x^2u_{N-1}(t,0) \\
\partial_x^2u_N(t,0)
\end{array}\right]_{N\times 1}=\left[\begin{array}{c}
0\\
0\\
\vdots\\
0  
\end{array}\right]_{2N\times 1} ; \quad 
\left[ C_4 \right ]_{2N \times N}\left[\begin{array}{c}
\partial_x^3 u_1(t,0)\\
\partial_x^3u_2(t,0)\\
\partial_x^3u_3(t,0) \\
\vdots\\
\partial_x^3u_{N-1}(t,0) \\
\partial_x^3u_N(t,0)
\end{array}\right]_{N\times 1}=\left[\begin{array}{c}
0\\
0\\
\vdots\\
0  
\end{array}\right]_{2N\times 1},
\end{equation}
where
\begin{equation*}
\left[ C_1 \right ]_{2N \times N}:=
  \begin{array}{c@{\!\!\!}l}
        \underbrace{ \left[\begin{array}{rrrrrr}
1          &-1          & 0         & \cdots  & 0          & 0    \\
0          & 1          & -1        &  \cdots & 0          & 0    \\
\vdots & \vdots & \vdots & \vdots  & \vdots & \vdots       \\
0          &  0        & 0          & \cdots  & 1          & -1     \\
\\
0          &  0        & 0          & \cdots  & 0          & 0     \\
\vdots & \vdots & \vdots & \vdots  & \vdots & \vdots       \\
0          &  0        & 0          & \cdots  & 0          & 0     
\end{array}\right]}_{\text{\tiny $N$ columns}}
&
 \begin{array}[c]{@{}l@{\,}l}
   \left. \begin{array}{c} \vphantom{0}  \\ \vphantom{\vdots} \\
   \\ \vphantom{0} \end{array} \right\} & \text{ \tiny $N-1$ rows} \\
   \\
\left. \begin{array}{c} \vphantom{0} \\ \vphantom{\vdots}
   \\ \vphantom{0}  \end{array} \right\} & \text{\tiny $N+1$ rows} 
\end{array}
\end{array}
\end{equation*}

\begin{equation*}
\left[ C_2 \right ]_{2N \times N}:=\begin{array}{c@{\!\!\!}l}
\underbrace{\left[\begin{array}{rrrrrr}
0          &0          & 0         & \cdots  & 0          & 0    \\
0          & 0          & 0        &  \cdots & 0          & 0    \\
\vdots & \vdots & \vdots & \vdots  & \vdots & \vdots       \\
0          &  0        & 0          & \cdots  & 0        & 0     \\
\\
1          &-1          & 0         & \cdots  & 0          & 0    \\
0          & 1          & -1        &  \cdots & 0          & 0    \\
\vdots & \vdots & \vdots & \vdots  & \vdots & \vdots       \\
0          &  0        & 0          & \cdots  & 1          & -1     \\
\\
0          &  0        & 0          & \cdots  & 0       & 0   \\
0          &  0        & 0          & \cdots  & 0         & 0     
\end{array}\right]}_{\text{\tiny $N$ columns}}
&
 \begin{array}[c]{@{}l@{\,}l}
   \left. \begin{array}{c} \vphantom{0}   \\ \vphantom{0} \\ 
   \\ \vphantom{0} \\  \end{array} \right\} & \text{\tiny $N-1$ rows} \\
   \\
   \\
\left. \begin{array}{c} \vphantom{0}  \\  \vphantom{0}  \\ 
   \\ \vphantom{0} \\ \end{array} \right\} & \text{\tiny $N-1$ rows} \\
   \\
   \left. \begin{array}{c} \vphantom{0} \\  \vphantom{0} \\ \end{array} \right\} & \text{\tiny $2$ rows} 
 \\
\end{array}
\end{array}
\end{equation*}
\begin{equation*}
\left[ C_3 \right ]_{2N \times N}:=\begin{array}{c@{\!\!\!}l}
\underbrace{\left[\begin{array}{rrrrrr}
0          &0          & 0         & \cdots  & 0          & 0    \\
0          & 0          & 0        &  \cdots & 0          & 0    \\
\vdots & \vdots & \vdots & \vdots  & \vdots & \vdots       \\
0          &  0        & 0          & \cdots  & 0          & 0     \\
\\
1          &  1        & 1         & \cdots  & 1       & 1   \\
0          &  0        & 0          & \cdots  & 0         & 0     
\end{array}\right]}_{\text{\tiny $N$ columns}}
&
 \begin{array}[c]{@{}l@{\,}l}
   \left. \begin{array}{c} \vphantom{0}   \\ \vphantom{0} \\\vphantom{0} \\
   \\ \vphantom{0} \\  \end{array} \right\} & \text{\tiny $2N-2$ rows} \\
   \\
 
\left. \begin{array}{c} \vphantom{0}  \\ 
   \\ \end{array} \right\} & \text{\tiny $2$ rows}
\end{array}
\end{array} 
\end{equation*}
and
\begin{equation*}
\left[ C_4 \right ]_{2N \times N}:=\begin{array}{c@{\!\!\!}l}
\underbrace{\left[\begin{array}{rrrrrr}
0          &0          & 0         & \cdots  & 0          & 0    \\
0          & 0          & 0        &  \cdots & 0          & 0    \\
0          & 0          & 0        &  \cdots & 0          & 0    \\
\vdots & \vdots & \vdots & \vdots  & \vdots & \vdots       \\
0          &  0        & 0          & \cdots  & 0         & 0     \\
\\
0          &  0        & 0         & \cdots  &0      & 0   \\
1          &  1        & 1         & \cdots  & 1       & 1   
\end{array}\right]}_{\text{\tiny $N$ columns}}
&
 \begin{array}[c]{@{}l@{\,}l}
   \left. \begin{array}{c} \vphantom{0}   \\ \vphantom{0} \\\vphantom{0} \\
   \\ \vphantom{0}    \\ \vphantom{0}\\  \end{array} \right\} & \text{\tiny $2N-2$ rows} \\
   \\
 
\left. \begin{array}{c} \vphantom{0}  \\ 
   \\ \end{array} \right\} & \text{\tiny $2$ rows}
\end{array}
\end{array}.
\end{equation*}
On the other hand, let be $\widetilde{u}_{j0}$ an extension for all line $\R$ of $u_{j0}$, satisfying
\begin{equation*}
\|\widetilde{u}_{j0}\|_{H^s(\R)}\lesssim\|u_{j0}\|_{H^s(\R^{+})}, \ j=1,2, ..., N,
\end{equation*}
respectively. Initially, we look for solutions in the form
\begin{align}\label{form}
&u_j(t,x)=\mathcal{L}^{\lambda_{j1}}\gamma_{j1}(t,x)+\mathcal{L}^{\lambda_{j2}}\gamma_{j2}(t,x)+F_j(t,x),\ j=1,2, ..., N.
\end{align}
Here, $\gamma_{ji}(\cdot)$,  $j=1,2, ..., N$, $i=1,2$, are unknown functions and
\begin{align*}
&F_j(t,x)=e^{it\partial_x^4}\widetilde{u}_{j0}+\mathcal{D}(\psi_T|u_j|^2u_j)(t,x),
\end{align*}
where $\mathcal{D}(w(t,x))=-i\int_0^t e^{i(t-t')\partial_x^4}w(t',x)dt'$. By using Lemma \ref{trace} we see that
\begin{equation}\label{tec1}
\begin{split}
u_j(t,0)=&\frac{M}{8\ }\left(\frac{e^{-i\frac{\pi}{8}(1+3\lambda_{j1})}+e^{-i\frac{\pi}{8}(1-5\lambda_{j1})}}{\sin(\frac{1-\lambda_{j1}}{4}\pi)}\right)\gamma_{j1}(t)+\frac{M}{8\ }\left(\frac{e^{-i\frac{\pi}{8}(1+3\lambda_{j2})}+e^{-i\frac{\pi}{8}(1-5\lambda_{j2})}}{\sin(\frac{1-\lambda_{j2}}{4}\pi)}\right)\gamma_{j2}(t)\\
&+F_j(t,0),\\
:=& a_{j1}\gamma_{j1}(t)+a_{j2}\gamma_{j2}(t)+F_j(t,0), \quad j=1,2, ..., N.
\end{split}
\end{equation}
Now, let us calculate the traces of first derivative functions. Thanks to \eqref{gamma}, \eqref{conv}, \eqref{forcing2} and Lemma \ref{trace}, we get that
\begin{equation}\label{tec2}
\begin{split}
\partial_xu_j(t,0)=&\frac{M}{8\ }\left(\frac{e^{-i\frac{\pi}{8}(-2+3\lambda_{j1})}+e^{-i\frac{\pi}{8}(6-5\lambda_{j1})}}{\sin(\frac{2-\lambda_{j1}}{4}\pi)}\right)\mathcal{I}_{-1/4}\gamma_{j1}(t)\\&+\frac{M}{8\ }\left(\frac{e^{-i\frac{\pi}{8}(-2+3\lambda_{j2})}+e^{-i\frac{\pi}{8}(6-5\lambda_{j2})}}{\sin(\frac{2-\lambda_{j2}}{4}\pi)}\right)\mathcal{I}_{-1/4}\gamma_{j2}(t)\\&+\partial_xF_j(t,0),\\
:=& b_{j1}\mathcal{I}_{-1/4}\gamma_{j1}(t)+b_{j2}\mathcal{I}_{-1/4}\gamma_{j2}(t)+\partial_x F_j(t,0), \quad j=1,2, ..., N.
\end{split}
\end{equation}
In the same way, we can have the traces of second and third derivatives functions, giving us the following
\begin{equation}\label{tec3}
\begin{split}
\partial_x^2u_j(t,0)=&
\frac{M}{8\ }\left(\frac{e^{-i\frac{\pi}{8}(-5+3\lambda_{j1})}+e^{-i\frac{\pi}{8}(11-5\lambda_{j1})}}{\sin(\frac{3-\lambda_{j1}}{4}\pi)}\right)\mathcal{I}_{-1/2}\gamma_{j1}(t)\\ 
&+\frac{M}{8\ }\left(\frac{e^{-i\frac{\pi}{8}(-5+3\lambda_{j2})}+e^{-i\frac{\pi}{8}(11-5\lambda_{j2})}}{\sin(\frac{3-\lambda_{j2}}{4}\pi)}\right)\mathcal{I}_{-1/2}\gamma_{j2}(t)\\
&+\partial_x^2F_j(t,0),\\
:=& c_{j1}\mathcal{I}_{-1/2}\gamma_{j1}(t)+c_{j2}\mathcal{I}_{-1/2}\gamma_{j2}(t)+\partial_x^2 F_j(t,0), \quad  j=1,2, ..., N
\end{split}
\end{equation}
and
\begin{equation}\label{tec4}
\begin{split}
\partial_x^3u_j(t,0)=& \frac{M}{8\ }\left(\frac{e^{-i\frac{\pi}{8}(-8+3\lambda_{j1})}+e^{-i\frac{\pi}{8}(16-5\lambda_{j1})}}{\sin(\frac{4-\lambda_{j1}}{4}\pi)}\right)\mathcal{I}_{-3/4}\gamma_{j1}(t)\\&+\frac{M}{8\ }\left(\frac{e^{-i\frac{\pi}{8}(-8+3\lambda_{j2})}+e^{-i\frac{\pi}{8}(16-5\lambda_{j2})}}{\sin(\frac{4-\lambda_{j2}}{4}\pi)}\right)\mathcal{I}_{-3/4}\gamma_{j2}(t)\\
&+\partial_x^2F_j(t,0),\\
:=&d_{j1}\mathcal{I}_{-3/4}\gamma_{j1}(t)+d_{j2}\mathcal{I}_{-3/4}\gamma_{j2}(t)+\partial_x^3 F_j(t,0), \quad  j=1,2, ..., N.
\end{split}
\end{equation}
Observe that Lemmas \ref{continuous} and \ref{trace} ensure these calculus are valid for $\mbox{Re} \lambda>-4$.
By substituting \eqref{tec1}, \eqref{tec2}, \eqref{tec3} and \eqref{tec4}   into \eqref{m1} and \eqref{m3}, yields that the functions $\gamma_{ji}$ and indexes $\lambda_{ji}$, for $j=1,2, ..., N$ and $i=1,2$, satisfy the following equalities:
\begin{equation*}
\left[ C_1 \right ]_{\scriptscriptstyle 2N \times N}\left[\begin{array}{ccccccc}
a_{11}&a_{12}& 0&0 & \cdots & 0 & 0\\
0& 0&a_{21}& a_{22} & \cdots & 0 & 0 \\
\vdots & \vdots & \vdots & \vdots  & \vdots & \vdots       \\
0& 0& 0& 0& \cdots & a_{N1}& a_{N2} 
\end{array}\right]_{\scriptscriptstyle N\times 2N}\left[\begin{array}{c}
\gamma_{11}(t)\\
\gamma_{12}(t)\\
\gamma_{21}(t)\\
\gamma_{22}(t)\\
\vdots \\ 
\gamma_{N1}(t)\\
\gamma_{N2}(t)
\end{array}\right]_{\scriptscriptstyle 2N\times 1} \\
=-\left[ C_1 \right ]_{\scriptscriptstyle 2N \times N}\left[\begin{array}{r}
F_1(t,0)\\
F_2(t,0)\\
\vdots \\
F_N(t,0)
\end{array}\right]_{\scriptscriptstyle N\times 1},
\end{equation*}
\begin{equation*}
\left[ C_2 \right ]_{\scriptscriptstyle 2N \times N}\left[\begin{array}{ccccccc}
b_{11}&b_{12}& 0&0 & \cdots & 0 & 0\\
0& 0&b_{21}& b_{22} & \cdots & 0 & 0 \\
\vdots & \vdots & \vdots & \vdots  & \vdots & \vdots       \\
0& 0& 0& 0& \cdots & b_{N1}& b_{N2} 
\end{array}\right]_{\scriptscriptstyle N\times 2N}\left[\begin{array}{c}
\gamma_{11}(t)\\
\gamma_{12}(t)\\
\gamma_{21}(t)\\
\gamma_{22}(t)\\
\vdots \\ 
\gamma_{N1}(t)\\
\gamma_{N2}(t)
\end{array}\right]_{\scriptscriptstyle 2N\times 1} \\
=-\left[ C_2 \right ]_{\scriptscriptstyle 2N \times N}\left[\begin{array}{r}
\partial_x\mathcal{I}_{\frac14}F_1(t,0)\\
\partial_x\mathcal{I}_{\frac14}F_2(t,0)\\
\vdots\\
\partial_x\mathcal{I}_{\frac14}F_N(t,0)
\end{array}\right]_{\scriptscriptstyle N\times 1},
\end{equation*}
\begin{equation*}
\left[ C_3 \right ]_{\scriptscriptstyle 2N \times N}\left[\begin{array}{ccccccc}
c_{11}&c_{12}& 0&0 & \cdots & 0 & 0\\
0& 0&c_{21}& c_{22} & \cdots & 0 & 0 \\
\vdots & \vdots & \vdots & \vdots  & \vdots & \vdots       \\
0& 0& 0& 0& \cdots & c_{N1}& c_{N2} 
\end{array}\right]_{\scriptscriptstyle N\times 2N}\left[\begin{array}{c}
\gamma_{11}(t)\\
\gamma_{12}(t)\\
\gamma_{21}(t)\\
\gamma_{22}(t)\\
\vdots \\ 
\gamma_{N1}(t)\\
\gamma_{N2}(t)
\end{array}\right]_{\scriptscriptstyle 2N\times 1} \\
=-\left[ C_3 \right ]_{\scriptscriptstyle 2N \times N}\left[\begin{array}{r}
\partial_x^2\mathcal{I}_{\frac12}F_1(t,0)\\
\partial_x^2\mathcal{I}_{\frac12}F_2(t,0)\\
\vdots \\
\partial_x^2\mathcal{I}_{\frac12}F_N(t,0)
\end{array}\right]_{\scriptscriptstyle N \times 1}
\end{equation*}
and
\begin{equation*}
\left[ C_4 \right ]_{\scriptscriptstyle 2N \times N}\left[\begin{array}{ccccccc}
d_{11}&d_{12}& 0&0 & \cdots & 0 & 0\\
0& 0&d_{21}& d_{22} & \cdots & 0 & 0 \\
\vdots & \vdots & \vdots & \vdots  & \vdots & \vdots       \\
0& 0& 0& 0& \cdots & d_{N1}& d_{N2} 
\end{array}\right]_{\scriptscriptstyle N\times 2N}\left[\begin{array}{c}
\gamma_{11}(t)\\
\gamma_{12}(t)\\
\gamma_{21}(t)\\
\gamma_{22}(t)\\
\vdots \\ 
\gamma_{N1}(t)\\
\gamma_{N2}(t)
\end{array}\right]_{\scriptscriptstyle 2N\times 1} \\
=-\left[ C_4 \right ]_{\scriptscriptstyle 2N \times N}\left[\begin{array}{c}
\partial_x^3\mathcal{I}_{\frac34}F_1(t,0)\\
\partial_x^3\mathcal{I}_{\frac34}F_2(t,0)\\
\vdots \\
\partial_x^3\mathcal{I}_{\frac34}F_N(t,0)
\end{array}\right]_{\scriptscriptstyle N \times 1}.
\end{equation*}
It follows that,
\begin{multline}\label{mother}
\left[\begin{array}{rrrrrrrrrrr}
	a_{11}& 	a_{12}& -a_{21}&-a_{22} & 0 & 0 & \cdots &0 &0 &0&0\\
	0&0&	a_{11}& 	a_{12}& -a_{21}&-a_{22} &  \cdots &0 &0 &0&0\\
	\vdots & 	\vdots & 	\vdots & 	\vdots & 	\vdots & 	\vdots & 	\vdots & 	\vdots & 	\vdots & 	\vdots & 	\vdots \\
	0&0&0&0&0&0&\cdots&	a_{(N-1)1}& 	a_{(N-1)2}& -a_{N1}&-a_{N2} \\
	\\
	b_{11}& 	b_{12}& -b_{21}&-b_{22} & 0 & 0 & \cdots &0 &0 &0&0\\
	0&0&	b_{11}& 	b_{12}& -b_{21}&-b_{22} &  \cdots &0 &0 &0&0\\
	\vdots & 	\vdots & 	\vdots & 	\vdots & 	\vdots & 	\vdots & 	\vdots & 	\vdots & 	\vdots & 	\vdots & 	\vdots \\
	0&0&0&0&0&0&\cdots&	b_{(N-1)1}& 	b_{(N-1)2}& -b_{N1}&-b_{N2} \\
	\\
	c_{11}&c_{12}&c_{21}&c_{22}&c_{31}&c_{32}&\cdots & \cdots&\cdots& c_{N1}&c_{N2} \\
	d_{11}&d_{12}&d_{21}&d_{22}&d_{31}&d_{32}&\cdots & \cdots &\cdots & d_{N1}&d_{N2} 
	\end{array}\right]_{\scriptscriptstyle 2N\times 2N}\left[\begin{array}{c}
\gamma_{11}(t)\\
\gamma_{12}(t)\\
\gamma_{21}(t)\\
\gamma_{22}(t)\\
\vdots \\ 
\gamma_{N1}(t)\\
\gamma_{N2}(t)
\end{array}\right]_{\scriptscriptstyle 2N\times 1}
\\
=-\left[\begin{array}{c}
F_1(t,0)-F_2(t,0)\\
\vdots \\
F_{N-1}(t,0)-F_N(t,0)\\
\\
\partial_x\mathcal{I}_{\frac14}F_1(t,0)-\partial_x\mathcal{I}_{\frac14}F_2(t,0)\\
\vdots \\
\partial_x\mathcal{I}_{\frac14}F_{N-1}(t,0)-\partial_x\mathcal{I}_{\frac14}F_{N}(t,0)\\
\\
\sum_{j=1}^N \partial_x^2\mathcal{I}_{\frac12}F_j(t,0)\\
\sum_{j=1}^N \partial_x^3\mathcal{I}_{\frac34}F_j(t,0)
\end{array}\right]_{\scriptscriptstyle 2N\times 1}.
\end{multline}
To simplify the notation,  let us denote the equality \eqref{mother} by
\begin{equation}\label{matrix}
\mathbf{M}(\lambda_{11},\lambda_{12}, \cdots, \lambda_{N1}, \lambda_{N2})\boldsymbol{\gamma}=\mathbf{F},
\end{equation}
where $\mathbf{M}(\lambda_{11},\lambda_{12}, \cdots,\lambda_{N1}, \lambda_{N2})$ is the first matrix that appears in the left hand side of \eqref{mother}, $\boldsymbol{\gamma}$ is the  matrix column given by vector $(\gamma_{11},\gamma_{12},\cdots,\gamma_{N1},\gamma_{N2})$ and $\mathbf{F}$ is the matrix in the right hand side of \eqref{mother}.

\subsection{Choosing the appropriate parameters and functions} In this second step, we need to choose the parameters $\lambda_{ji}$ and the functions $\gamma_{ji}$, with $j=1,2, ..., N$, $i=1,2$, in such a way that we can be able to write the solution $u_j(t,x)$, in a integral form.

To do this, let us start by using the hypothesis of Lemma \ref{operator}. We need, firstly, to fix parameters $\lambda_{ji}$, such that
\begin{equation}\label{lam}
\max\left\{\frac{2s-7}{2},-1\right\}<\lambda_{ji}(s)<\min\left\{s+\frac12,\frac12\right\}, \quad j=1,2, ..., N, i=1,2.
\end{equation}
With this restriction in hand  we choose the parameters $\lambda_{ji}$ as follows
\begin{equation}\label{choices}
\lambda_{11}=\lambda_{21}=\cdots=\lambda_{N1}=-\frac12\  \text{and}\  \lambda_{12}=\lambda_{22}=\cdots=\lambda_{N2}=\frac14,
\end{equation} 
then, we have the equation
\begin{equation}\label{matrix2}
\mathbf{M}\left(-\frac12,\frac14,\cdots, -\frac12, \frac14\right)\boldsymbol{\gamma}=\mathbf{F}.
\end{equation}
The following lemma gives us  that $\mathbf{M}\left(-\frac12,\frac14,	\cdots, -\frac12, \frac14\right)$ is invertible.

\begin{lemma}\label{tec}
	The determinant of matrix $\mathbf{M}\left(-\frac12,\frac14, \cdots, -\frac12, \frac14\right)$ is nonzero.
\end{lemma}

We will prove  Lemma \ref{tec}  on the next section. Thus, this good choices of the parameters satisfying \eqref{lam} together with this lemma ensures that $\mathbf{M}$ is invertible and, consequently, the following holds
\begin{equation}\label{gamma1}
\boldsymbol{\gamma}=\mathbf{M}^{-1}\left(-\frac12,\frac14,	\cdots, -\frac12, \frac14\right)\mathbf{F}.
\end{equation}
We empathize that $\gamma_{ji}$ depends on $F_1$ and $F_2$, which depend on the unknown functions $u_1$ and $u_2$. Thus, by substituting \eqref{gamma1} into \eqref{form}, we get $u_j(t,x)$ in the integral form
\begin{equation}\label{form1}
u_j(t,x)=\mathcal{L}^{-\frac12}\gamma_{j1}(t,x)+\mathcal{L}^{\frac14}\gamma_{j2}(t,x)+F_j(t,x),
\quad j=1,2, ..., N.
\end{equation}

\subsection{Defining the truncated integral operator and functional space}\label{truncated}
Using the previous subsection, we have the solution of the Cauchy problem \eqref{grapha} with Type $\mathcal{A}$ boundary condition \eqref{boundary1} in the integral form \eqref{form1}. In order to use the Fourier restriction method, the third step is to define a truncated version for the integral form \eqref{form1}. 

Pick $s\in[0,1/2)$, we fix the parameters $\lambda_{ji}$ as in \eqref{choices}  and define $$\boldsymbol{\gamma}=(\gamma_{11},\gamma_{12},\gamma_{21},\gamma_{22}, 	\cdots, \gamma_{N1},\gamma_{N2} )$$ by \eqref{gamma1}. Consider $b=b(s)<\frac12$  and that the estimates given in Lemmas \ref{grupo}, \ref{duhamel}, \ref{operator} and \ref{trilinear} are valid. Now, define the operator $\Lambda$ by
\begin{equation*}
\Lambda= \left(\Lambda_1,\Lambda_2, \cdots, \Lambda_N 	\right)
\end{equation*}
where
\begin{align*}
&\Lambda_j u(t,x)=\psi(t)\mathcal{L}^{-\frac12}\gamma_{j1}(t,x)+\psi(t)\mathcal{L}^{\frac14}\gamma_{j2}(t,x)+F_j(t,x), \quad j=1,2,...,N.
\end{align*}
Here,
\begin{align*}
&F_j(t,x)=\psi(t)(e^{it\partial_x^4}\widetilde{u}_{j0}+\lambda\mathcal{D}(\psi_T|u_j|^2u_j)(t,x)), \quad j=1,2,...,N,
\end{align*}
with $$\mathcal{D}(w(t,x))=-i\int_0^t e^{i(t-t')\partial_x^4}w(x,t')dt'.$$
We consider $\Lambda$ defined on the Banach space $Z(s,b)=\bigoplus_{j =1}^NZ_j(s,b)$ by
\begin{multline*}
Z_j(s,b)=\left\lbrace w\in C(\R_t;H^s(\R_x)) \cap C(\R_x;H^{\frac{2s+3}{8}}(\R_t))\cap X^{s,b}; \right.\\
 \left. w_x\in C(\R_x;H^{\frac{2s+1}{8}}(\R_t)),w_{xx}\in C(\R_x;H^{\frac{2s-1}{8}}(\R_t)), w_{xxx}\in C(\R_x;H^{\frac{2s-3}{8}}(\R_t)) \right\rbrace, 
\end{multline*}
for $j=1,2,..., N$, with norm
\begin{equation}\label{norma}
\|(u_1,u_2,\cdots,u_N)\|_{Z(s,b)}=\|u_1\|_{Z_1(s,b)}+\|u_2\|_{Z_2(s,b)}+\cdots+\|u_N\|_{Z_N(s,b)}.
\end{equation}
Each norm of $\|u\|_{Z_j(s,b)}$ on \eqref{norma} is defined by
\begin{multline*}
\|u\|_{Z_j(s,b)}=\|u\|_{C\left(\R_t;H^s(\R_x)\right)}+\|u\|_{C\left(\R_x;H^{\frac{2s+3}{8}}(\R_t)\right)}+\|u\|_{ X^{s,b}}\\
+\|u_x\|_{C\left(\R_x;H^{\frac{2s+1}{8}}(\R_t)\right)}+\|u_{xx}\|_{C\left(\R_x;H^{\frac{2s-1}{8}}(\R_t)\right)}+\|u_{xxx}\|_{C\left(\R_x;H^{\frac{2s-3}{8}}(\R_t)\right)},
\end{multline*}
for $j=1,2,..., N$.

\subsection{Proving that the functions $\mathcal{L}^{-\frac12}\gamma_{j1}$ and $\mathcal{L}^{\frac14}\gamma_{j2}$, for $j=1,2,..., N$, are well-defined}\label{welldefined}

%Now, we will prove that the functions $\mathcal{L}^{\lambda_{11}}\gamma_1(t,x)$ and $ \mathcal{L}^{\lambda_{12}}\gamma_2(t,x)$ are well defined, for any function $u=(u_1,u_2)\in Z$ . 
Indeed, by using Lemma \ref{operator} it suffices to show that these functions are in the closure of the spaces $C_0^{\infty}(\R^+)$. By using expression \eqref{gamma1}, we see that $\gamma_{ji}$ ($j=1,2,...,N$ and $ i=1,2$) are linear combinations of the functions 
\begin{equation*}
\begin{array}{c}
F_1(t,0)-F_2(t,0), \quad  F_2(t,0)-F_3(t,0), \quad \cdots ,  \quad F_{N-1}(t,0)-F_N(t,0),\\
\\
\partial_x\mathcal{I}_{\frac14}F_1(t,0)-\partial_x\mathcal{I}_{\frac14}F_2(t,0), \quad \partial_x\mathcal{I}_{\frac14}F_2(t,0)-\partial_x\mathcal{I}_{\frac14}F_3(t,0), \quad  \cdots, \quad \partial_x\mathcal{I}_{\frac14}F_{N-1}(t,0)-\partial_x\mathcal{I}_{\frac14}F_{N}(t,0),\\
\\
\partial_x^2\mathcal{I}_{\frac12}F_1(t,0)+\partial_x^2\mathcal{I}_{\frac12}F_2(t,0)+\cdots+\partial_x^2\mathcal{I}_{\frac12}F_N(t,0),\\
\\
 \partial_x^3\mathcal{I}_{\frac34}F_1(t,0)+ \partial_x^3\mathcal{I}_{\frac34}F_2(t,0)+\cdots+ \partial_x^3\mathcal{I}_{\frac34}F_N(t,0).
 \end{array}
\end{equation*}
Thus, we need to show that the functions $F_j(t,0),\  \partial_x\mathcal{I}_{\frac{1}{3}}F_j(t,0),\  \partial_x^2\mathcal{I}_{\frac{2}{3}}F_j(t,0)$ are in appropriate spaces. By Lemmas \ref{grupo}, \ref{operator}, \ref{duhamel} and \ref{trilinear} we obtain
\begin{equation}\label{trace1}
\|F_j(t,0)\|_{H^{\frac{2s+3}{8}}(\R^+)}\leq c( \|u_{j0}\|_{H^s(\R^+)}+\|u_j\|_{X^{s,b}}^3).
\end{equation}
If $0\leq s<\frac{1}{2}$ we have that $\frac38 <\frac{2s+3}{8}<\frac12$, then $H^{\frac{2s+3}{8}}(\R^+)=H_0^{\frac{2s+3}{8}}(\R^+)$. It follows that $F_j(t,0)\in H_0^{\frac{2s+3}{8}}(\R^+)$ for $0\leq s<\frac{1}{2}.$
Again by using Lemmas  \ref{grupo}, \ref{operator}, \ref{duhamel} and \ref{trilinear} we get
\begin{equation*}
\|\partial_x F_j(t,0)\|_{H^{\frac{2s+1}{8}}(\R^+)}\leq c( \|u_{j0}\|_{H^s(\R^+)}+\|u_j\|_{X^{s,b}}^3). 
\end{equation*}
Since $0\leq s<\frac12$ we have $\frac18\leq \frac{2s+1}{8}<\frac14$, then  the functions $\partial_xF_j(t,0)\in H_0^{\frac{2s+1}{8}}(\R^+)$. Then, thanks to Lemma \ref{lio}, we have that
\begin{equation*}
\|\partial_x \mathcal{I}_{\frac14}F_j(t,0)\|_{H_0^{\frac{2s+3}{8}}(\R^+)}\leq c( \|u_{j0}\|_{H^s(\R^+)}+\|u_j\|_{X^{s,b}}^3).
\end{equation*}
Therefore, this yields that
\begin{equation}\label{trace2}
\partial_x\mathcal{I}_{\frac{1}{4}}F_j(t,0)-\partial_x\mathcal{I}_{\frac{1}{4}}F_{j+1}(t,0)\in H_0^{\frac{2s+3}{8}}(\R^+), \quad j=1,2,..., N.
\end{equation}
In a similar way, we can obtain
\begin{equation*}
\begin{split}
\|\partial_x^2\mathcal{I}_{\frac{1}{2}}F_j(t,0)\|_{H_0^{\frac{2s+3}{8}}(\R^+)}\lesssim \|u_{j0}\|_{H^s(\R^+)}+\|u_j\|_{X^{s,b}}^3,\\
\|\partial_x^3\mathcal{I}_{\frac{3}{4}}F_i(t,0)\|_{H_0^{\frac{2s+3}{8}}(\R^+)}\lesssim \|u_{j0}\|_{H^s(\R^+)}+\|u_j\|_{X^{s,b}}^3.
\end{split}
\end{equation*}
It follows that
\begin{equation}\label{trace3}
\begin{split}
\partial_x^2\mathcal{I}_{\frac{1}{2}}F_1(t,0)+\partial_x^2\mathcal{I}_{\frac{1}{2}}F_2(t,0)+\cdots+ \partial_x^2\mathcal{I}_{\frac{1}{2}}F_N(t,0) \in H_0^{\frac{2s+3}{8}}(\R^+),\\
\partial_x^3\mathcal{I}_{\frac{3}{4}}F_1(t,0)+\partial_x^3\mathcal{I}_{\frac{3}{4}}F_2(t,0)+\cdots + \partial_x^3\mathcal{I}_{\frac{3}{4}}F_N(t,0)\in H_0^{\frac{2s+3}{8}}(\R^+).
\end{split}
\end{equation}
Thus, \eqref{trace1}, \eqref{trace2} and \eqref{trace3} imply that the functions  $\mathcal{L}^{-\frac12}\gamma_{j1}$ and $\mathcal{L}^{\frac14}\gamma_{j2}$, for $j=1,2,..., N$, are well-defined.

\subsection{Showing that $\Lambda$ is a contraction in a ball of $Z$}\label{contraction}
Lemmas  \ref{lio}, \ref{grupo}, \ref{operator}, \ref{duhamel} and \ref{trilinear} guarantee that
\begin{multline*}
\|\Lambda(u_{1},\cdots, u_N)-\Lambda(v_{1},\cdots, v_{N})\|_{Z_{s,b}} \\
\leq T^{\epsilon} c\left(\|(u_1,\cdots,u_N)\|_Z^2+\|(v_1,\cdots, v_N)\|_{Z}^2 \right)\|(u_1,\cdots, u_N)-(v_1,\cdots, v_N)\|_{Z}
\end{multline*}
and 
\begin{equation*}
\begin{split}
\|\Lambda(u_1,\cdots, u_N)\|_{Z_{s,b}}\leq c& \left( \|u_{01}\|_{H^s(\R^+)}+\cdots+\|u_{0N}\|_{H^s(\R^+)}\right)+T^{\epsilon}(\|u_1\|_{X^{s,b}}^3+\cdots+\|u_N\|_{X^{s,b}}^3),
\end{split}
\end{equation*}
for $\epsilon$ adequately small.

Consider in $Z$ the ball defined by 
$$
B=\{ (u_1,\cdots, u_N)\in Z^{s,b}; \|(u_1,\cdots, u_N)\|_{Z^{s,b}}\leq M\},
$$
where $$M=2 c\left(\|u_{01}\|_{H^s(\R^+)}+\cdots+\|u_{0N}\|_{H^s(\R^+)}\right).$$ Lastly, choosing $T=T(M)$ sufficiently small, such that 
$$\|\Lambda (u_1,\cdots, u_N)\|_{Z^{s,b}}\leq M$$
and $$\|\Lambda (u_1,\cdots, u_N) - \Lambda (v_1,\cdots, v_N)\|_{Z^{s,b}}\leq \frac12\norm{(u_1,\cdots, u_N)-(v_1,\cdots, v_N)}_{Z^{s,b}},$$ it follows that $\Lambda$ is a contraction map on $B$ and has a fixed point $(\widetilde{u}_1, \cdots, \widetilde{u}_N)$.  Therefore, the restriction
\begin{equation*}
(u_1,\cdots, u_N)=(\widetilde{u}_1\big|_{\R^-\times (0,T)}, \cdots, \widetilde{u}_N\big|_{\R^+\times (0,T)})
\end{equation*}
solves the Cauchy problem \eqref{grapha} with Type $\mathcal A$ vertex boundary condition \eqref{boundary1}, in the sense of distributions. Thus, Theorem \ref{theorem1} is a consequence of the above steps, described in the previous subsections, finalizing so the proof. \qed

\section{Proof of lemma \ref{tec}}\label{sec5}
First, we will prove the case $N=2$. The vertex conditions \eqref{boundary1}, for this case, is given by
 \begin{equation*}
\begin{cases}
\partial_x^k u_1(t,0)=\partial_x^k u_2(t,0),\ k=0,1 & t\in(0,T)\\
\sum_{j=1}^2\partial_x^ku_j(t,0)=0,\ k=2,3\ & t\in(0,T).
\end{cases}
\end{equation*}
In this way, we consider the vertex conditions as the following matrices
\begin{equation}\label{m111}
\left[\begin{array}{cc}
1&-1\\
0& 0 \\
0& 0 \\
0& 0 \\
\end{array}\right]\left[\begin{array}{r}
u_1(t,0)\\
u_2(t,0)
\end{array}\right]=0, \quad \left[\begin{array}{cc}
0&0\\
1& -1 \\
0& 0 \\
0& 0 \\
\end{array}\right]\left[\begin{array}{r}
\partial_xu_1(t,0)\\
\partial_x u_2(t,0)
\end{array}\right]=0
\end{equation}
and
\begin{equation}\label{m311}
\left[\begin{array}{cc}
0&0\\
0& 0 \\
1& 1 \\
0& 0 \\
\end{array}\right]\left[\begin{array}{r}
\partial_x^2u_1(t,0)\\
\partial_x^2u_2(t,0)
\end{array}\right]=0, \quad  \left[\begin{array}{cc}
0&0\\
0& 0 \\
0& 0 \\
1& 1 \\
\end{array}\right]\left[\begin{array}{r}
\partial_x^3u_1(t,0)\\
\partial_x^3u_2(t,0)
\end{array}\right]=0.
\end{equation}

By substituting, for $N=2$, \eqref{tec1}, \eqref{tec2},  \eqref{tec3} and \eqref{tec4} into \eqref{m111} and \eqref{m311}, yields that the functions $\gamma_{ji}$ and indexes $\lambda_{ji}$, for $j=1,2$ and $i=1,2$, satisfy the equality of matrices:
\begin{equation*}
\begin{split}
\left[\begin{array}{cc}
1&-1\\
0& 0 \\
0&  0\\
0& 0 \\
\end{array}\right]\left[\begin{array}{cccc}
a_{11}&a_{12}& 0&0\\
0& 0&a_{21}& a_{22} \\
\end{array}\right]\left[\begin{array}{r}
\gamma_{11}(t)\\
\gamma_{12}(t)\\
\gamma_{21}(t)\\
\gamma_{22}(t)\\
\end{array}\right]&=-\left[\begin{array}{cc}
1&-1\\
0& 0\\
0& 0\\
0& 0 \\
\end{array}\right]\left[\begin{array}{r}
F_1(t,0)\\
F_2(t,0)\\
\end{array}\right],
\\
\left[\begin{array}{cc}
0&0\\
1& -1 \\
0& 0 \\
0& 0 \\
\end{array}\right]\left[\begin{array}{cccc}
b_{11}&b_{12}& 0&0\\
0& 0&b_{21}& b_{22} \\
\end{array}\right]\left[\begin{array}{r}
\gamma_{11}(t)\\
\gamma_{12}(t)\\
\gamma_{21}(t)\\
\gamma_{22}(t)\\
\end{array}\right]&=-\left[\begin{array}{cc}
0&0\\
1& -1 \\
0& 0 \\
0& 0\\
\end{array}\right]\left[\begin{array}{r}
\partial_x\mathcal{I}_{\frac14}F_1(t,0)\\
\partial_x\mathcal{I}_{\frac14}F_2(t,0)\\
\end{array}\right],
\\
\left[\begin{array}{cc}
0&0\\
0& 0 \\
1& 1 \\
0& 0 \\
\end{array}\right]\left[\begin{array}{cccc}
c_{11}&c_{12}& 0&0\\
0& 0&c_{21}& c_{22} \\
\end{array}\right]\left[\begin{array}{r}
\gamma_{11}(t)\\
\gamma_{12}(t)\\
\gamma_{21}(t)\\
\gamma_{22}(t)\\
\end{array}\right]&=-\left[\begin{array}{cc}
0&0 \\
0& 0 \\
1& 1 \\
0& 0\\
\end{array}\right]\left[\begin{array}{r}
\partial_x^2\mathcal{I}_{\frac12}F_1(t,0)\\
\partial_x^2\mathcal{I}_{\frac12}F_2(t,0)\\
\end{array}\right]
\end{split}
\end{equation*}
and
\begin{equation*}
\left[\begin{array}{cc}
0&0\\
0& 0 \\
0& 0 \\
1& 1 \\
\end{array}\right]\left[\begin{array}{cccc}
d_{11}&d_{12}& 0&0\\
0& 0&d_{21}& d_{22} \\
\end{array}\right]\left[\begin{array}{r}
\gamma_{11}(t)\\
\gamma_{12}(t)\\
\gamma_{21}(t)\\
\gamma_{22}(t)\\
\end{array}\right]=-\left[\begin{array}{cccc}
0&0\\
0& 0 \\
0& 0 \\
1& 1 \\
\end{array}\right]\left[\begin{array}{r}
\partial_x^3\mathcal{I}_{\frac34}F_1(t,0)\\
\partial_x^3\mathcal{I}_{\frac34}F_2(t,0)\\
\end{array}\right].
\end{equation*}
Putting all matrices together, we have that,
\begin{equation*}
\begin{split}
&\left[\begin{array}{cccc}
	a_{11}&
	a_{12}& -a_{21}&-a_{22}\\
\\	b_{11}&
	b_{12}& -b_{21}&-b_{22}\\
\\
	c_{11}&c_{12}&c_{21}&c_{22}\\
	\\
		d_{11}&d_{12}&d_{21}&d_{22}\\
	\end{array}\right]\left[\begin{array}{r}
	\gamma_{11}\\
\\	
	\gamma_{12}\\
	\\
	\gamma_{21}\\
	\\
	\gamma_{22}
	\end{array}\right]
=-\left[\begin{array}{c}
F_1(t,0)-F_2(t,0)\\
\\
\partial_x\mathcal{I}_{\frac14}F_1(t,0)-\partial_x\mathcal{I}_{\frac14}F_2(t,0)\\
\\
\partial_x^2\mathcal{I}_{\frac12}F_1(t,0)+\partial_x^2\mathcal{I}_{\frac12}F_2(t,0)\\
\\
\partial_x^3\mathcal{I}_{\frac34}F_1(t,0)+\partial_x^3\mathcal{I}_{\frac34}F_2(t,0)
\end{array}\right].
\end{split}
\end{equation*}
In the case $N=2$, the matrix $M$, given by \eqref{mother}, can be read as follows
\begin{equation}\label{mothercaseN=2}
M=
\left[\begin{array}{cccc}
	a_{11}&
	a_{12}& -a_{21}&-a_{22}\\
\\	b_{11}&
	b_{12}& -b_{21}&-b_{22}\\
\\
	c_{11}&c_{12}&c_{21}&c_{22}\\
	\\
		d_{11}&d_{12}&d_{21}&d_{22}\\
	\end{array}\right],
\end{equation}
where $a_{ij}, b_{ij}, c_{ij}$ and $d_{ij}$ are given by \eqref{tec1}, \eqref{tec2}, \eqref{tec3} and \eqref{tec4}, respectively.

\begin{claim}\label{claim1}
$M$ has determinant different of zero with appropriate choice of $\lambda_{ji}$, $j=1,2$ and $i=1,2$. 
\end{claim}

In fact, firstly noting that $$\sin\left(\frac{2-a}{4}\pi\right)=\cos\left(\frac{a\pi}{4}\right)$$ and it is easy to see that
\begin{equation}\label{coefi}
\begin{cases}
a_{ji}=\frac{M e^{-\frac{i\pi}{8}}}{8}\left( \frac{e^{-\frac{3i\pi\lambda_{ji}}{8}}+e^{\frac{5i\pi\lambda_{ji}}{8}}}{  \sin\left(\frac{(1-\lambda_{ji})\pi}{4}\right)	}\right), &     b_{ji}=\frac{M e^{\frac{2i\pi}{8}}}{8}\left( \frac{e^{-\frac{3i\pi\lambda_{ji}}{8}}-e^{\frac{5i\pi\lambda_{ji}}{8}}}{ \cos\left(\frac{\lambda_{ji}\pi}{4}\right) 	}\right), \\
\\
c_{ji}=\frac{ M e^{\frac{5i\pi}{8}}}{8}\left( \frac{e^{-\frac{3i\pi\lambda_{ji}}{8}}+e^{\frac{5i\pi\lambda_{ji}}{8}}}{ \sin\left(\frac{(3-\lambda_{ji})\pi}{4}\right)	}\right), & d_{ji}=\frac{M e^{\frac{8i\pi}{8}}}{8}\left( \frac{e^{-\frac{3i\pi\lambda_{ji}}{8}}-e^{\frac{5i\pi\lambda_{ji}}{8}}}{ \sin\left(\frac{\lambda_{ji}\pi}{4}\right) 	}\right).
\end{cases} \quad j=1,2, \quad i=1,2
\end{equation}
Then, the determinant of $M$ can be write as 

\begin{equation*}
| M|=\frac{M^4 e^{\frac{14i\pi}{8}}}{8^4}\left|\begin{array}{cccc} 
 \frac{e^{-\frac{3i\pi\lambda_{11}}{8}}+e^{\frac{5i\pi\lambda_{11}}{8}}}{\sin\left(\frac{(1-\lambda_{11})\pi}{4}\right)	} &  \frac{e^{-\frac{3i\pi\lambda_{12}}{8}}+e^{\frac{5i\pi\lambda_{12}}{8}}}{\sin\left(\frac{(1-\lambda_{12})\pi}{4}\right)	} & - \frac{e^{-\frac{3i\pi\lambda_{21}}{8}}+e^{\frac{5i\pi\lambda_{21}}{8}}}{\sin\left(\frac{(1-\lambda_{21})\pi}{4}\right)	} &  -\frac{e^{-\frac{3i\pi\lambda_{22}}{8}}+e^{\frac{5i\pi\lambda_{22}}{8}}}{\sin\left(\frac{(1-\lambda_{22})\pi}{4}\right)} \\
 \\
 \frac{e^{-\frac{3i\pi\lambda_{11}}{8}}-e^{\frac{5i\pi\lambda_{11}}{8}}}{ \cos\left(\frac{\lambda_{11}\pi}{4}\right) 	}& \frac{e^{-\frac{3i\pi\lambda_{12}}{8}}-e^{\frac{5i\pi\lambda_{12}}{8}}}{ \cos\left(\frac{\lambda_{12}\pi}{4}\right) 	} & - \frac{e^{-\frac{3i\pi\lambda_{21}}{8}}-e^{\frac{5i\pi\lambda_{21}}{8}}}{ \cos\left(\frac{\lambda_{21}\pi}{4}\right) 	} & - \frac{e^{-\frac{3i\pi\lambda_{22}}{8}}-e^{\frac{5i\pi\lambda_{22}}{8}}}{ \cos\left(\frac{\lambda_{22}\pi}{4}\right) 	}\\
\\
 \frac{e^{-\frac{3i\pi\lambda_{11}}{8}}+e^{\frac{5i\pi\lambda_{11}}{8}}}{\sin\left(\frac{(3-\lambda_{11})\pi}{4}\right)	} &  \frac{e^{-\frac{3i\pi\lambda_{12}}{8}}+e^{\frac{5i\pi\lambda_{12}}{8}}}{\sin\left(\frac{(3-\lambda_{12})\pi}{4}\right)	} & \frac{e^{-\frac{3i\pi\lambda_{21}}{8}}+e^{\frac{5i\pi\lambda_{21}}{8}}}{\sin\left(\frac{(3-\lambda_{21})\pi}{4}\right)	} &  \frac{e^{-\frac{3i\pi\lambda_{22}}{8}}+e^{\frac{5i\pi\lambda_{22}}{8}}}{\sin\left(\frac{(3-\lambda_{22})\pi}{4}\right)} \\
\\
 \frac{e^{-\frac{3i\pi\lambda_{11}}{8}}-e^{\frac{5i\pi\lambda_{11}}{8}}}{ \sin\left(\frac{\lambda_{11}\pi}{4}\right) 	} &  \frac{e^{-\frac{3i\pi\lambda_{12}}{8}}-e^{\frac{5i\pi\lambda_{12}}{8}}}{ \sin\left(\frac{\lambda_{12}\pi}{4}\right) 	}&  \frac{e^{-\frac{3i\pi\lambda_{21}}{8}}-e^{\frac{5i\pi\lambda_{21}}{8}}}{ \sin\left(\frac{\lambda_{21}\pi}{4}\right) 	} &  \frac{e^{-\frac{3i\pi\lambda_{22}}{8}}-e^{\frac{5i\pi\lambda_{22}}{8}}}{ \sin\left(\frac{\lambda_{22}\pi}{4}\right) 	}
\end{array}\right|. 
\end{equation*}
By using the identity $$	\frac{e^{ia}-e^{ib}}{e^{ia}+e^{ib}}=i\  \text{tan}\left(\frac{a-b}{2}\right ),$$ we have that
\begin{equation*}
\begin{split}
	\left| M\left(\lambda_{11},\lambda_{12},\lambda_{21},\lambda_{22}\right) \right|=\frac{M^4 e^{\frac{14i\pi}{8}}}{8^4}& \times \left\lbrace (e^{-\frac{3\lambda_{11}\pi}{8}i}+e^{\frac{5\lambda_{11}\pi}{8}i})(e^{-\frac{3\lambda_{12}\pi}{8}i}+e^{\frac{5\lambda_{12}\pi}{8}i}) \right\rbrace \times\\
&\left\lbrace (e^{-\frac{3\lambda_{21}\pi}{8}i}+e^{\frac{5\lambda_{21}\pi}{8}i})   (e^{-\frac{3\lambda_{22}\pi}{8}i}+e^{\frac{5\lambda_{22}\pi}{8}i})\right\rbrace \left| A\right|,
\end{split}
\end{equation*}
	where $A$ is the matrix
	\begin{equation*}
A=	\left[\begin{array}{cccc} \frac{1}{\text{sin}\left(\frac{1-\lambda_{11}}{4}\pi\right)}& \frac{1}{\text{sin}\left(\frac{1-\lambda_{12}}{4}\pi\right)} & -\frac{1}{\text{sin}\left(\frac{1-\lambda_{21}}{4}\pi\right)} & -\frac{1}{\text{sin}\left(\frac{1-\lambda_{22}}{4}\pi\right)} \\
\\
 -\frac{i \text{tan}(\frac{\lambda_{11}}{2}\pi)}{\text{cos}(\frac{\lambda_{11}}{4}\pi)} & -\frac{i \text{tan}(\frac{\lambda_{12}}{2}\pi)}{\text{cos}(\frac{\lambda_{12}}{4}\pi)} &  \frac{i \text{tan}(\frac{\lambda_{21}}{2}\pi)}{\text{cos}(\frac{\lambda_{21}}{4}\pi)}&  \frac{i \text{tan}(\frac{\lambda_{22}}{2}\pi)}{\text{cos}(\frac{\lambda_{22}}{4}\pi)}   \\
\\
  \frac{1}{\text{sin}\left(\frac{3-\lambda_{11}}{4}\pi\right)}& \frac{1}{\text{sin}\left(\frac{3-\lambda_{12}}{4}\pi\right)} & \frac{1}{\text{sin}\left(3-\frac{\lambda_{21}}{4}\pi\right)} & \frac{1}{\text{sin}\left(\frac{3-\lambda_{22}}{4}\pi\right)} \\
  \\ 
  -\frac{i \text{tan}(\frac{\lambda_{11}}{2}\pi)}{\text{sin}(\frac{\lambda_{11}}{4}\pi)} & -\frac{i \text{tan}(\frac{\lambda_{12}}{2}\pi)}{\text{sin}(\frac{\lambda_{12}}{4}\pi)} & - \frac{i \text{tan}(\frac{\lambda_{21}}{2}\pi)}{\text{sin}(\frac{\lambda_{21}}{4}\pi)}&-  \frac{i \text{tan}(\frac{\lambda_{22}}{2}\pi)}{\text{sin}(\frac{\lambda_{22}}{4}\pi)}\end{array}\right].
	\end{equation*}
	
	By choosing $\lambda_{11}=\lambda_{21}$ and $\lambda_{12}=\lambda_{22}$ we have that the constant that appears before of the matrix $A$ takes the form:
	\begin{equation*}
	\frac{2M^4e^{\frac{14}{8}\pi i}(e^{-i\lambda_{11}\pi}+1)^2(e^{-i\lambda_{12}{\pi}}+1)^2e^{\frac{5}{4}\pi(\lambda_{11}+\lambda_{12})i}}{8^4}.
	\end{equation*}
	Note that this number is zero only in the case $ \lambda_{11} = 2k+1$ and $\lambda_{12}=2l+1$  for $k,l \in \mathbb{Z}$.

	Let us denote the entries of the matrix $A$ as follows:
	\begin{equation}\label{model1}
	A=\left[\begin{array}{cccc} a& n & -a & -n \\ f & g & -f&  -g\\ c& e & c & e   \\ d &m &  d& m\end{array}\right].
	\end{equation}
Thus, its determinant is given by 
	\begin{equation}\label{det1}
	\text{det} A=4(de-cm)(nf-ag).
	\end{equation}

	In particular for $ \lambda_{11}=\lambda_{21}=-\frac12$ and $\lambda_{12}=\lambda_{22}=\frac14$ matrix $A$ can be seen as follows,
	\begin{equation*}
A'=	\left[\begin{array}{cccc} \sqrt 2\sqrt{2-\sqrt{2}}& \sqrt 2 \sqrt{4 - 2 \sqrt{2} + \sqrt{2 (10 - 7 \sqrt{2})}} & -\sqrt 2\sqrt{2-\sqrt{2}} & -\sqrt 2\sqrt{4 - 2 \sqrt{2} + \sqrt{2 (10 - 7 \sqrt{2})}} \\ \frac{2i}{\sqrt{2+\sqrt{2}}} & \frac{i\text tan(\pi/8)}{\text{cos}(\pi/16)} & -\frac{2i}{\sqrt{2+\sqrt{2}}}&  -\frac{i\text tan(\pi/8)}{\text{cos}(\pi/16)}\\ \sqrt 2\sqrt{2+\sqrt{2}}& \sqrt 2 \sqrt{4 - 2 \sqrt{2} - \sqrt{2 (10 - 7 \sqrt{2})}}  & \sqrt 2 \sqrt{2+\sqrt{2}} & \sqrt 2 \sqrt{4 - 2 \sqrt{2} - \sqrt{2 (10 - 7 \sqrt{2})}}    \\ -i\sqrt{2(2+\sqrt2)} &- \frac{i\text{tan}(\pi/8)}{\text{sin}(\frac{\pi}{16})}&   -i\sqrt{2(2+\sqrt2)} &  - \frac{i\text{tan}(\pi/8)}{\text{sin}(\frac{\pi}{16})}\end{array}\right].
	\end{equation*}
	
	By using determinant properties the determinant of $A'$ is equivalent of the determinant of the following matrix:
	\begin{equation*}
	-(4i)(\sqrt{2}i)\left[\begin{array}{cccc} \sqrt{2-\sqrt{2}}& \sqrt{4 - 2 \sqrt{2} + \sqrt{2 (10 - 7 \sqrt{2})}} & -\sqrt{2-\sqrt{2}} & -\sqrt{4 - 2 \sqrt{2} + \sqrt{2 (10 - 7 \sqrt{2})}} \\ \frac{1}{\sqrt{2+\sqrt{2}}} & \frac{\text tan(\pi/8)}{2\text{cos}(\pi/16)} & -\frac{1}{\sqrt{2+\sqrt{2}}}&  -\frac{\text tan(\pi/8)}{2\text{cos}(\pi/16)}\\ \sqrt{2+\sqrt{2}}&  \sqrt{4 - 2 \sqrt{2} -\sqrt{2 (10 - 7 \sqrt{2})}}  & \sqrt{2+\sqrt{2}} &  \sqrt{4 - 2 \sqrt{2} - \sqrt{2 (10 - 7 \sqrt{2})}}    \\ \sqrt{2+\sqrt2} & \frac{\text{tan}(\pi/8)}{\sqrt2\text{sin}(\frac{\pi}{16})}&   \sqrt{2+\sqrt2} &   \frac{\text{tan}(\pi/8)}{\sqrt2\text{sin}(\frac{\pi}{16})}\end{array}\right].
	\end{equation*}
	Therefore, we can rewrite  matrix \eqref{model1} as follows
	\begin{equation*}
	\left[\begin{array}{cccc} a& n & -a & -n \\ \frac{1}{c} & g & -\frac{1}{c}&  -g\\ c& e & c & e   \\ c &m &  c&  m\end{array}\right]
	\end{equation*}
	and its determinant is given by 
	\begin{equation*}
	4(e-m)(n-acg).
	\end{equation*}
	We only need to check that $e-m\neq 0$ and $n-acg\neq0$. 
	An calculation proves that 
	$$e-m\sim-0,6508$$   
	and
	$$n-acg\sim0,9741. $$
	Thus, we have that
\begin{equation}\label{new6}
det(A')=-4(2i)(\sqrt2i)(e-m)(n-acg)\sim-7,1722
\end{equation}
%	
%	\begin{equation*}
%	\begin{split}
%	&\left(\sqrt{2(2+\sqrt{2-\sqrt{2}})}+\frac{2(2 \sqrt{2}-3)}{\sqrt{2-\sqrt{2+\sqrt{2}}}}\right)\left(\sqrt{4-2 \sqrt{2}+\sqrt{2(10-7 \sqrt{2})}}-\sqrt{\frac{(2-\sqrt{2})(2+\sqrt{2})}{2+\sqrt{2+\sqrt{2}}}}\right)\\
%	&\sim0.81241202104726071787940516703961221867593424854969541558 \ldots
%	\end{split}
%	\end{equation*}
that is,  the determinant of matrix $\mathbf{M}\left(-\frac12,\frac14, -\frac12, \frac14\right)$ given by \eqref{mothercaseN=2} is nonzero, proving the Claim \ref{claim1} and Lemma \ref{tec}, for the case $N=2$.

\vspace{0.2cm}

For a better understanding of the reader, before to do the general case, we will present briefly also the proof of  Lemma \ref{tec} considering $N = 3$. For instance, vertex conditions  \eqref{boundary1}, in this case, are given in the matrices form as follows:
\begin{equation*}
\left[\begin{array}{ccc}
1&-1&0\\
0& 1&-1\\
0& 0 &0\\
0& 0 &0\\
0& 0 &0\\
0& 0 &0
\end{array}\right]\left[\begin{array}{r}
u_1(t,0)\\
u_2(t,0)\\
u_3(t,0)
\end{array}\right]= 0, \quad \quad
\left[\begin{array}{ccc}
0&0&0\\
1& -1 &0\\
0& 1&-1 \\
0& 0 &0\\
0& 0 &0\\
0& 0 &0\\
\end{array}\right]\left[\begin{array}{r}
\partial_xu_1(t,0)\\
\partial_x u_2(t,0)\\
\partial_x u_3(t,0)
\end{array}\right]=0,
\end{equation*}
and
\begin{equation*}
\left[\begin{array}{ccc}
0&0&0\\
0& 0 &0 \\
1& 1 &1\\
0& 0 & 0\\
0& 0 & 0\\
0& 0 & 0\\
\end{array}\right]\left[\begin{array}{r}
\partial_x^2u_1(t,0)\\
\partial_x^2u_2(t,0)\\
\partial_x^2u_3(t,0)
\end{array}\right]=0,
\quad \quad
\left[\begin{array}{ccc}
0&0&0\\
0& 0&0 \\
0& 0&0 \\
1& 1 &1\\
0& 0 & 0\\
0& 0 & 0\\
\end{array}\right]\left[\begin{array}{r}
\partial_x^3u_1(t,0)\\
\partial_x^3u_2(t,0)\\
\partial_x^3u_3(t,0)
\end{array}\right]=0.
\end{equation*}
Thus, combining the above matrices and the integral form of solution \eqref{form}, as in the case $N=2$, we obtain
\begin{equation*}
\left[\begin{array}{cccccc}
	a_{11}&
    a_{12}& -a_{21}&-a_{22}&0&0\\
\\	0&
	0& a_{21}&a_{22}&-a_{31}&-a_{32}\\
\\	b_{11}&
	b_{12}& -b_{21}&-b_{22}&0&0\\
\\	0&
	0& b_{21}&b_{22}&-b_{31}&-b_{32}\\
\\
	c_{11}&c_{12}&c_{21}&c_{22}&c_{31}&c_{32}\\
	\\
		d_{11}&d_{12}&d_{21}&d_{22}&d_{31}&d_{32}\\
	\end{array}\right]\left[\begin{array}{r}
	\gamma_{11}\\
\\	
	\gamma_{12}\\
	\\
	\gamma_{21}\\
	\\
	\gamma_{22}\\
	\\
	\gamma_{31}\\
	\\
	\gamma_{32}
	\end{array}\right]
=-\left[\begin{array}{c}
F_1(t,0)-F_2(t,0)\\
\\
F_2(t,0)-F_3(t,0)\\
\\
\partial_x\mathcal{I}_{\frac14}F_1(t,0)-\partial_x\mathcal{I}_{\frac14}F_2(t,0)\\
\\
\partial_x\mathcal{I}_{\frac14}F_2(t,0)-\partial_x\mathcal{I}_{\frac14}F_3(t,0)\\
\\
\partial_x^2\mathcal{I}_{\frac12}F_1(t,0)+\partial_x^2\mathcal{I}_{\frac12}F_2(t,0)+\partial_x^2\mathcal{I}_{\frac12}F_3(t,0)\\
\\
\partial_x^3\mathcal{I}_{\frac34}F_1(t,0)+\partial_x^3\mathcal{I}_{\frac34}F_2(t,0)+\partial_x^3\mathcal{I}_{\frac34}F_3(t,0)
\end{array}\right].
\end{equation*}

Let us consider $M$ the following matrix
\begin{equation}\label{mothercaseN=3}
M=\left[\begin{array}{cccccc}
	a_{11}&
    a_{12}& -a_{21}&-a_{22}&0&0\\
\\	0&
	0& a_{21}&a_{22}&-a_{31}&-a_{32}\\
\\	b_{11}&
	b_{12}& -b_{21}&-b_{22}&0&0\\
\\	0&
	0& b_{21}&b_{22}&-b_{31}&-b_{32}\\
\\
	c_{11}&c_{12}&c_{21}&c_{22}&c_{31}&c_{32}\\
	\\
		d_{11}&d_{12}&d_{21}&d_{22}&d_{31}&d_{32}\\
	\end{array}\right].
\end{equation}

\begin{claim}\label{claim2}
$M$ has determinant different of zero with appropriate choice of $\lambda_{ji}$, $j=1,2,3$ and $i=1,2$. 
\end{claim}

Indeed, similarly as in the case $N=2$ and by using the identities \eqref{coefi}, yields that
\begin{equation*}
\begin{cases}
a_{ji}=\frac{M e^{-\frac{i\pi}{8}}}{8}\left( \frac{e^{-\frac{3i\pi\lambda_{ji}}{8}}+e^{\frac{5i\pi\lambda_{ji}}{8}}}{  \sin\left(\frac{(1-\lambda_{ji})\pi}{4}\right)	}\right), &     b_{ji}=\frac{M e^{\frac{2i\pi}{8}}}{8}\left( \frac{e^{-\frac{3i\pi\lambda_{ji}}{8}}-e^{\frac{5i\pi\lambda_{ji}}{8}}}{ \cos\left(\frac{\lambda_{ji}\pi}{4}\right) 	}\right), \\
\\
c_{ji}=\frac{ M e^{\frac{5i\pi}{8}}}{8}\left( \frac{e^{-\frac{3i\pi\lambda_{ji}}{8}}+e^{\frac{5i\pi\lambda_{ji}}{8}}}{ \sin\left(\frac{(3-\lambda_{ji})\pi}{4}\right)	}\right), & d_{ji}=\frac{M e^{\frac{8i\pi}{8}}}{8}\left( \frac{e^{-\frac{3i\pi\lambda_{ji}}{8}}-e^{\frac{5i\pi\lambda_{ji}}{8}}}{ \sin\left(\frac{\lambda_{ji}\pi}{4}\right) 	}\right).
\end{cases} \quad j=1,2,3, \quad i=1,2.
\end{equation*}
By determinant properties, we can get the determinant of $M$ as
\begin{align*}
|M|&=\left(\frac{Me^{-\frac{ i\pi}{8}}}{8} \right)^2\left(\frac{Me^{\frac{2 i\pi}{8}}}{8} \right)^2\left(\frac{Me^{\frac{ 5i\pi}{8}}}{8} \right)\left(\frac{Me^{\frac{ 8i\pi}{8}}}{8} \right)\prod_{i=1,2, j=1,2,3}  \left( e^{-\frac{3i\pi\lambda_{ji}}{8}}+e^{\frac{5i\pi\lambda_{ji}}{8}}\right)|A| \\
&= \frac{M^6e^{\frac{ 15i\pi}{8}}}{8^6}  \prod_{i=1,2, j=1,2,3}  \left( e^{-\frac{3i\pi\lambda_{ji}}{8}}+e^{\frac{5i\pi\lambda_{ji}}{8}}\right)|A|,
\end{align*}
	where $A$ is the following matrix
	\begin{equation*}
A'=	\left[\begin{array}{cccccc} 
\frac{1}{ \sin \left( \frac{1-\lambda_{11}}{4} \pi \right) } & \frac{1}{ \sin \left( \frac{1-\lambda_{12}}{4} \pi \right) }&- \frac{1}{ \sin \left( \frac{1-\lambda_{21}}{4} \pi \right) } & - \frac{1}{ \sin \left( \frac{1-\lambda_{22}}{4} \pi \right) } & 0 & 0\\
\\
0& 0& \frac{1}{ \sin \left( \frac{1-\lambda_{21}}{4} \pi \right) } & \frac{1}{ \sin \left( \frac{1-\lambda_{22}}{4} \pi \right) }&- \frac{1}{ \sin \left( \frac{1-\lambda_{31}}{4} \pi \right) } & - \frac{1}{ \sin \left( \frac{1-\lambda_{32}}{4} \pi \right) } \\
\\
-\frac{i \text{tan}(\frac{\lambda_{11}}{2}\pi)}{\text{cos}(\frac{\lambda_{11}}{4}\pi)} & -\frac{i \text{tan}(\frac{\lambda_{12}}{2}\pi)}{\text{cos}(\frac{\lambda_{12}}{4}\pi)} &  \frac{i \text{tan}(\frac{\lambda_{21}}{2}\pi)}{\text{cos}(\frac{\lambda_{21}}{4}\pi)}&  \frac{i \text{tan}(\frac{\lambda_{22}}{2}\pi)}{\text{cos}(\frac{\lambda_{22}}{4}\pi)}  & 0 & 0 \\
\\
0&0&-\frac{i \text{tan}(\frac{\lambda_{21}}{2}\pi)}{\text{cos}(\frac{\lambda_{21}}{4}\pi)} & -\frac{i \text{tan}(\frac{\lambda_{22}}{2}\pi)}{\text{cos}(\frac{\lambda_{22}}{4}\pi)} &  \frac{i \text{tan}(\frac{\lambda_{31}}{2}\pi)}{\text{cos}(\frac{\lambda_{31}}{4}\pi)}&  \frac{i \text{tan}(\frac{\lambda_{32}}{2}\pi)}{\text{cos}(\frac{\lambda_{32}}{4}\pi)}  \\
\\
\frac{1}{ \sin \left( \frac{3-\lambda_{11}}{4} \pi \right) } & \frac{1}{ \sin \left( \frac{3-\lambda_{12}}{4} \pi \right) }& \frac{1}{ \sin \left( \frac{3-\lambda_{21}}{4} \pi \right) } &  \frac{1}{ \sin \left( \frac{3-\lambda_{22}}{4} \pi \right) } &  \frac{1}{ \sin \left( \frac{3-\lambda_{31}}{4} \pi \right) } &  \frac{1}{ \sin \left( \frac{3-\lambda_{32}}{4} \pi \right) }\\
  \\ 
  -\frac{i \text{tan}(\frac{\lambda_{11}}{2}\pi)}{\text{sin}(\frac{\lambda_{11}}{4}\pi)} & -\frac{i \text{tan}(\frac{\lambda_{12}}{2}\pi)}{\text{sin}(\frac{\lambda_{12}}{4}\pi)} & - \frac{i \text{tan}(\frac{\lambda_{21}}{2}\pi)}{\text{sin}(\frac{\lambda_{21}}{4}\pi)}&-  \frac{i \text{tan}(\frac{\lambda_{22}}{2}\pi)}{\text{sin}(\frac{\lambda_{22}}{4}\pi)} &-  \frac{i \text{tan}(\frac{\lambda_{31}}{2}\pi)}{\text{sin}(\frac{\lambda_{31}}{4}\pi)} & -  \frac{i \text{tan}(\frac{\lambda_{32}}{2}\pi)}{\text{sin}(\frac{\lambda_{32}}{4}\pi)}
  \end{array}\right].
	\end{equation*}

	By choosing $\lambda_{11}=\lambda_{21}=\lambda_{31}$ and $\lambda_{12}=\lambda_{22}=\lambda_{32}$, we have that the constant that appears before of the matrix $A$ takes the form:
	\begin{equation*}
	\frac{M^6e^{\frac{15 i\pi}{8}}}{8^6} (e^{-i\pi \lambda_{11}}+1)^3(e^{-i\pi \lambda_{12}}+1)^3e^{\frac{15}{4}\pi(\lambda_{11}+\lambda_{12})i}.
	\end{equation*}
	Note that this number is zero only in the case $\lambda_{11} = 2n+1$ and $ \lambda_{12}=2m+1$  for $n,m \in \mathbb{Z}$. 	Let us  rewrite the entries of matrix $A$ as follows:
	\begin{equation*}\label{model1_1} 
	A=
	\left[\begin{array}{cccccc} a& n & -a & -n & 0 & 0
	 \\
	 0 & 0 & a& n & -a & -n 
	 \\
	  f & g & -f&  -g & 0 & 0
	   \\
	0 & 0 &  f & g & -f&  -g 
	  \\ 
	  c& e & c & e  &c & e
	   \\
	   d &m &  d&  m & d & m
	  \end{array}\right].
	\end{equation*}
Thus, its determinant is given by 
	\begin{equation*}
|A'|=9(de-cm)(ag-nf)^2.
	\end{equation*}
	
Finally, considering  $\lambda_{11}=\lambda_{21}=\lambda_{31}=-\frac12$ and $\lambda_{12}=\lambda_{22}=\lambda_{32}=\frac14$, thanks to the case $N=2$, we have that  $(ag-fn)\neq 0$ and $(de-cm)\neq0$, thus $|A|\neq 0$. Claim \ref{claim2} is thus proved and Lemma \ref{tec1} is achieved, when $N=3$.

\vglue 0.2cm

 Let us now deal with the general situation, that is, when $N>3$. Consider $$\mathbf{ M}=\mathbf{M}(\lambda_{11},\lambda_{12}, \cdots, \lambda_{N1}, \lambda_{N2})$$ defined by \eqref{mother}, namely,
 \begin{equation*}
\mathbf{M}_{\scriptscriptstyle 2N\times 2N}=\begin{array}{c@{\!\!\!}l}
\underbrace{\left[\begin{array}{rrrrrrrrrrr}
	a_{11}& 	a_{12}& -a_{21}&-a_{22} & 0 & 0 & \cdots &0 &0 &0&0\\
	0&0&	a_{11}& 	a_{12}& -a_{21}&-a_{22} &  \cdots &0 &0 &0&0\\
	\vdots & 	\vdots & 	\vdots & 	\vdots & 	\vdots & 	\vdots & 	\vdots & 	\vdots & 	\vdots & 	\vdots & 	\vdots \\
	0&0&0&0&0&0&\cdots&	a_{(N-1)1}& 	a_{(N-1)2}& -a_{N1}&-a_{N2} \\
	\\
	b_{11}& 	b_{12}& -b_{21}&-b_{22} & 0 & 0 & \cdots &0 &0 &0&0\\
	0&0&	b_{11}& 	b_{12}& -b_{21}&-b_{22} &  \cdots &0 &0 &0&0\\
	\vdots & 	\vdots & 	\vdots & 	\vdots & 	\vdots & 	\vdots & 	\vdots & 	\vdots & 	\vdots & 	\vdots & 	\vdots \\
	0&0&0&0&0&0&\cdots&	b_{(N-1)1}& 	b_{(N-1)2}& -b_{N1}&-b_{N2} \\
	\\
	c_{11}&c_{12}&c_{21}&c_{22}&c_{31}&c_{32}&\cdots & \cdots&\cdots& c_{N1}&c_{N2} \\
	d_{11}&d_{12}&d_{21}&d_{22}&d_{31}&d_{32}&\cdots & \cdots &\cdots & d_{N1}&d_{N2} 
	\end{array}\right]}_{\text{\tiny $2N$ columns}}
&
 \begin{array}[c]{@{}l@{\,}l}
   \left. \begin{array}{c} \vphantom{0}   \\ \vphantom{0} \\ 
   \\ \vphantom{0}   \\ \vphantom{0}  \end{array} \right\} & \text{\tiny $N-1$ rows} \\
   \vphantom{}  \\
\left. \begin{array}{c} \vphantom{0}  \\  \vphantom{0}  \\ 
   \\ \vphantom{0} \\ \end{array} \right\} & \text{\tiny $N-1$ rows} \\
   \\
   \left. \begin{array}{c} \vphantom{0} \\  \vphantom{0} \\ \end{array} \right\} & \text{\tiny $2$ rows} 
 \\
\end{array}
\end{array},
\end{equation*}
where $a_{ij}, b_{ij}, c_{ij}$ and $d_{ij}$ are given by \eqref{tec1}, \eqref{tec2}, \eqref{tec3} and \eqref{tec4}, respectively. As we noted in the cases $N=2,3$,  this function of $\lambda_{ji}$ can be take the form 
\begin{equation*}
\begin{cases}
a_{ji}=\frac{M e^{-\frac{i\pi}{8}}}{8}\left( \frac{e^{-\frac{3i\pi\lambda_{ji}}{8}}+e^{\frac{5i\pi\lambda_{ji}}{8}}}{  \sin\left(\frac{(1-\lambda_{ji})\pi}{4}\right)	}\right), &     b_{ji}=\frac{M e^{\frac{2i\pi}{8}}}{8}\left( \frac{e^{-\frac{3i\pi\lambda_{ji}}{8}}-e^{\frac{5i\pi\lambda_{ji}}{8}}}{ \cos\left(\frac{\lambda_{ji}\pi}{4}\right) 	}\right), \\
\\
c_{ji}=\frac{ M e^{\frac{5i\pi}{8}}}{8}\left( \frac{e^{-\frac{3i\pi\lambda_{ji}}{8}}+e^{\frac{5i\pi\lambda_{ji}}{8}}}{ \sin\left(\frac{(3-\lambda_{ji})\pi}{4}\right)	}\right), & d_{ji}=\frac{M e^{\frac{8i\pi}{8}}}{8}\left( \frac{e^{-\frac{3i\pi\lambda_{ji}}{8}}-e^{\frac{5i\pi\lambda_{ji}}{8}}}{ \sin\left(\frac{\lambda_{ji}\pi}{4}\right) 	}\right).
\end{cases} \quad j=1,2,..., N, \quad i=1,2.
\end{equation*}
 Thus, by using the determinant properties, we have that 
\begin{align*}
|\mathbf{ M}|&=\left(\frac{Me^{-\frac{ i\pi}{8}}}{8} \right)^{N-1}\left(\frac{Me^{\frac{2 i\pi}{8}}}{8} \right)^{N-1}\left(\frac{Me^{\frac{ 5i\pi}{8}}}{8} \right)\left(\frac{Me^{\frac{ 8i\pi}{8}}}{8} \right)\prod_{i=1,2, j=1,\cdots,N}  \left( e^{-\frac{3i\pi\lambda_{ji}}{8}}+e^{\frac{5i\pi\lambda_{ji}}{8}}\right)|\mathbf{ M'}| \\
&= \frac{M^{2N}e^{\frac{ (12+N)i\pi}{8}}}{8^{2N}}  \prod_{i=1,2, j=1,\cdots,N}  \left( e^{-\frac{3i\pi\lambda_{ji}}{8}}+e^{\frac{5i\pi\lambda_{ji}}{8}}\right)|\mathbf{ M'}|,
\end{align*}
where $M'$ is a matrix, depending only of $\lambda_{ji}$, given by
  \begin{equation*}
\mathbf{M'}=\left[\begin{array}{rrrrrrrrrrr}
	\bar a_{11}& 	\bar a_{12}& - \bar a_{21}&- \bar a_{22} & 0 & 0 & \cdots &0 &0 &0&0\\
	0&0&	\bar a_{11}& 	\bar a_{12}& - \bar a_{21}&- \bar a_{22} &  \cdots &0 &0 &0&0\\
	\vdots & 	\vdots & 	\vdots & 	\vdots & 	\vdots & 	\vdots & 	\vdots & 	\vdots & 	\vdots & 	\vdots & 	\vdots \\
	0&0&0&0&0&0&\cdots&	\bar  a_{(N-1)1}& 	\bar a_{(N-1)2}& -\bar a_{N1}&-\bar a_{N2} \\
	\\
	\bar b_{11}& 	\bar b_{12}& -\bar b_{21}&-\bar b_{22} & 0 & 0 & \cdots &0 &0 &0&0\\
	0&0&	\bar b_{11}& 	\bar b_{12}& -\bar b_{21}&-\bar b_{22} &  \cdots &0 &0 &0&0\\
	\vdots & 	\vdots & 	\vdots & 	\vdots & 	\vdots & 	\vdots & 	\vdots & 	\vdots & 	\vdots & 	\vdots & 	\vdots \\
	0&0&0&0&0&0&\cdots&	\bar b_{(N-1)1}& 	\bar b_{(N-1)2}& -\bar b_{N1}&-\bar b_{N2} \\
	\\
	\bar c_{11}&\bar c_{12}&\bar c_{21}&\bar c_{22}&\bar c_{31}&\bar c_{32}&\cdots & \cdots&\cdots&\bar  c_{N1}&\bar c_{N2} \\
	\bar d_{11}&\bar  d_{12}&\bar d_{21}&\bar d_{22}&\bar  d_{31}&\bar d_{32}&\cdots & \cdots &\cdots &\bar  d_{N1}&\bar d_{N2} 
	\end{array}\right]_{\scriptscriptstyle 2N\times 2N}.
\end{equation*}
Here, the coefficients of matrix $\mathbf{M'}$ are given by
\begin{equation}\label{new5}
\begin{cases}
\bar a_{ji}=  \dfrac{1}{  \sin\left(\frac{(1-\lambda_{ji})\pi}{4}\right)	}, &    \bar b_{ji}= -\dfrac{i\tan \left( \frac{\lambda_{ji\pi}}{2}\right)}{ \cos\left(\frac{\lambda_{ji}\pi}{4}\right) 	}, \\
\\
\bar c_{ji}=\dfrac{1}{ \sin\left(\frac{(3-\lambda_{ji})\pi}{4}\right)}, & \bar d_{ji}= -\dfrac{i\tan \left( \frac{\lambda_{ji\pi}}{2}\right)}{ \sin\left(\frac{\lambda_{ji}\pi}{4}\right) 	}.
\end{cases} \quad j=1,2,..., N, \quad i=1,2.
\end{equation}
	By choosing $\lambda_{11}=\lambda_{21} =\cdots =\lambda_{N1}$ and $\lambda_{12}=\lambda_{22}=\cdots=\lambda_{N2}$, we have that the constant that appears before of the matrix $\mathbf{M'}$ takes the form:
	\begin{equation*}
	\frac{M^{2N}e^{\frac{(12+N) i\pi}{8}}}{8^{2N}} (e^{-i\pi \lambda_{11}}+1)^N(e^{-i\pi \lambda_{12}}+1)^Ne^{\frac{5N}{4}\pi(\lambda_{11}+\lambda_{12})i}.
	\end{equation*}
	Note that this number is zero only in the case $\lambda_{11} = 2n+1$ and $ \lambda_{12}=2m+1$  for $n,m \in \mathbb{Z}$. 	Let us denote the entries of matrix $\mathbf{ M}'$ as follows:
	\begin{equation*}
\mathbf{M'}=\left[\begin{array}{rrrrrrrrrrr}
	a& 	n& - a&- n & 0 & 0 & \cdots &0 &0 &0&0\\
	0&0&	a& 	n& - a&- n &  \cdots &0 &0 &0&0\\
	\vdots & 	\vdots & 	\vdots & 	\vdots & 	\vdots & 	\vdots & 	\vdots & 	\vdots & 	\vdots & 	\vdots & 	\vdots \\
	0&0&0&0&0&0&\cdots&	a& n& -a&-n \\
	\\
	f& 	g& -f&-g & 0 & 0 & \cdots &0 &0 &0&0\\
	0&0&	f& 	g& -f&-g &  \cdots &0 &0 &0&0\\
	\vdots & 	\vdots & 	\vdots & 	\vdots & 	\vdots & 	\vdots & 	\vdots & 	\vdots & 	\vdots & 	\vdots & 	\vdots \\
	0&0&0&0&0&0&\cdots&	f& 	g& -f&-g \\
	\\
	c&e&c&e&c&e&\dots & \cdots&\cdots&  c& e \\
	d&m&d&m&d&m& \dots & \cdots&\cdots&  d& m 
	\end{array}\right]_{\scriptscriptstyle 2N\times 2N}.
\end{equation*}
Moreover, by using the determinant properties, it yields that 
\begin{equation*}
\left| \mathbf{M'} \right|=\left|\begin{array}{rrrrrrrrrrr}
	a& 	n& - a&- n & 0 & 0 & \cdots &0 &0 &0&0\\
	f& 	g& -f&-g & 0 & 0 & \cdots &0 &0 &0&0\\
	0&0&	a& 	n& - a&- n &  \cdots &0 &0 &0&0\\
	0&0&	f& 	g& -f&-g &  \cdots &0 &0 &0&0\\
	\vdots & 	\vdots & 	\vdots & 	\vdots & 	\vdots & 	\vdots & 	\vdots & 	\vdots & 	\vdots & 	\vdots & 	\vdots \\
	0&0&0&0&0&0&\cdots&	a& n& -a&-n \\
	0&0&0&0&0&0&\cdots&	f& 	g& -f&-g \\
	\\
	c&e&c&e&c&e&\dots & \cdots&\cdots&  c& e \\
	d&m&d&m&d&m& \dots & \cdots&\cdots&  d& m 
	\end{array}\right|_{\scriptscriptstyle 2N\times 2N}.
\end{equation*}
Considering the matrix 
\begin{equation*}
A=\left[ \begin{array}{c c}
a & n \\
f & g
\end{array}\right] \quad \text{and} \quad 
B=\left[ \begin{array}{c c}
c & e \\
d & m
\end{array}\right],
\end{equation*}
the determinant of $\mathbf{M}'$ can be write as a block matrices, namely,
\begin{equation}\label{deter1}
\left| \mathbf{M'} \right|=\left| \begin{array}{cccccc}
A_{\scriptscriptstyle 2\times 2} &-A_{\scriptscriptstyle 2\times 2}  & 0_{\scriptscriptstyle 2\times 2} & \cdots &0_{\scriptscriptstyle 2\times 2} & 0_{\scriptscriptstyle 2\times 2} \\
 0_{\scriptscriptstyle 2\times 2} & A_{\scriptscriptstyle 2\times 2} &-A_{\scriptscriptstyle 2\times 2}  &  \cdots & 0_{\scriptscriptstyle 2\times 2} & 0_{\scriptscriptstyle 2\times 2} 	\\
 \vdots &  \vdots &   \vdots &   \vdots &   \vdots &   \vdots \\
 0_{\scriptscriptstyle 2\times 2} & 0_{\scriptscriptstyle 2\times 2} & 0_{\scriptscriptstyle 2\times 2} & \cdots & A_{\scriptscriptstyle 2\times 2} &-A_{\scriptscriptstyle 2\times 2} \\
 \\
 B_{\scriptscriptstyle 2\times 2} &  B_{\scriptscriptstyle 2\times 2}  &  B_{\scriptscriptstyle 2\times 2} & \cdots & B_{\scriptscriptstyle 2\times 2}  & B_{\scriptscriptstyle 2\times 2} 
\end{array}\right|_{\scriptscriptstyle 2N \times 2N}.
\end{equation}
From now on, we denote $0_{\scriptscriptstyle n\times n} $ and $I_{\scriptscriptstyle n\times n} $ the null and identity matrices, respectively.

Let us  introduce the properties of determinants that helped us to prove  Lemma \ref{tec1} in general form. Consider a block matrix $N$ of size $(n+m)\times (n+m)$ of the form
\begin{equation*}
N=\left[ \begin{array}{c c}
C & D \\
F & G
\end{array}\right],
\end{equation*}
where $C,D,F$ and $G$ are of size $n\times n$, $n\times m$, $m\times n$ and $m\times m$, respectively. If $G$ is invertible, then
\begin{equation}\label{blockdeter}
\det N=\det (C-DG^{-1}F)\det(G).
\end{equation}
In fact, this property follows immediately from the following identity 
\begin{equation*}
\left[ \begin{array}{c c}
C & D \\
F & G
\end{array}\right] \left[ \begin{array}{c c}
I & 0 \\
-G^{-1}F & I
\end{array}\right]=\left[ \begin{array}{c c}
C-DG^{-1}F & D \\
0 &  G
\end{array}\right].
\end{equation*}
Finally, recall that the determinant of a block triangular matrix is the product of the determinants of its diagonal blocks.

With these two properties in hand, define $C,D,F$ and $G$, respectively, by
\begin{equation*}
C=\left[ \begin{array}{ccccc}
A_{\scriptscriptstyle 2\times 2} &-A_{\scriptscriptstyle 2\times 2}  & 0_{\scriptscriptstyle 2\times 2} & \cdots &0_{\scriptscriptstyle 2\times 2}  \\
 0_{\scriptscriptstyle 2\times 2} & A_{\scriptscriptstyle 2\times 2} &-A_{\scriptscriptstyle 2\times 2}  &  \cdots & 0_{\scriptscriptstyle 2\times 2} 	\\
 \vdots &  \vdots &   \vdots &   \vdots &   \vdots  \\
 0_{\scriptscriptstyle 2\times 2} & 0_{\scriptscriptstyle 2\times 2} & 0_{\scriptscriptstyle 2\times 2} & \cdots & A_{\scriptscriptstyle 2\times 2} 
\end{array}\right]_{\scriptscriptstyle 2(N-2) \times 2(N-2)}, \quad 
D= \left[ \begin{array}{c}
0_{\scriptscriptstyle 2\times 2}   \\
 0_{\scriptscriptstyle 2\times 2} 	\\
 \vdots \\
 -A_{\scriptscriptstyle 2\times 2} 
\end{array}\right]_{\scriptscriptstyle  2(N-2)\times 2}
\end{equation*}
and
\begin{equation*}
F=\left[\begin{array}{ccccc}
 B_{\scriptscriptstyle 2\times 2} &  B_{\scriptscriptstyle 2\times 2}  &  B_{\scriptscriptstyle 2\times 2} & \cdots & B_{\scriptscriptstyle 2\times 2} 
\end{array}\right]_{\scriptscriptstyle 2 \times 2(N-2)}, \quad G=B_{\scriptscriptstyle 2\times 2}.
\end{equation*}
Thanks to the case $N=2$, we already know that  
\begin{equation}\label{det3}
\det G =  \det  B_{\scriptscriptstyle 2 \times 2}=cm - de \neq 0,
\end{equation}
which implies that $G$ is invertible. Thus, the determinant \eqref{deter1} takes the form 
\begin{equation*}
\left| \mathbf{M'} \right|=\left| \begin{array}{cc}
C_{\scriptscriptstyle 2(N-2)\times 2(N-2)} & D_{\scriptscriptstyle 2 \times 2(N-2)} \\
F_{\scriptscriptstyle 2 \times 2(N-2)} &  B_{\scriptscriptstyle 2 \times 2}
\end{array}\right|_{\scriptscriptstyle 2N \times 2N}
\end{equation*}
and by using the property \eqref{blockdeter}, it yields that 
\begin{equation}\label{det2}
\det \mathbf{M'} =\det \left( C_{\scriptscriptstyle 2(N-2)\times 2(N-2)}  - D_{\scriptscriptstyle 2 \times 2(N-2)} \left(B_{\scriptscriptstyle 2 \times 2}\right)^{-1}F_{\scriptscriptstyle 2 \times 2(N-2)}\right) \det  B_{\scriptscriptstyle 2 \times 2}.
\end{equation}

\begin{claim} \label{claim3}
$\mathbf{M}'$ has determinant different of zero with appropriate choice of $\lambda_{ji}$, $j=1,2,..., N$ and $i=1,2$.
\end{claim}
From \eqref{det3} is enough to prove that $$\det \left( C_{\scriptscriptstyle 2(N-2)\times 2(N-2)}  - D_{\scriptscriptstyle 2 \times 2(N-2)} \left(B_{\scriptscriptstyle 2 \times 2}\right)^{-1}F_{\scriptscriptstyle 2 \times 2(N-2)}\right)$$ is  nonzero. In order to analyze the above determinant, note that 
\begin{equation*}
\begin{split}
\left(B_{\scriptscriptstyle 2 \times 2}\right)^{-1}F_{\scriptscriptstyle 2 \times 2(N-2)}&=\left(B_{\scriptscriptstyle 2 \times 2}\right)^{-1}\left[\begin{array}{ccccc}
 B_{\scriptscriptstyle 2\times 2} &  B_{\scriptscriptstyle 2\times 2}  &  B_{\scriptscriptstyle 2\times 2} & \cdots & B_{\scriptscriptstyle 2\times 2}\end{array}\right]_{\scriptscriptstyle 2 \times 2(N-2)} \\
 & = \left[\begin{array}{ccccc}
 I_{\scriptscriptstyle 2\times 2} &  I_{\scriptscriptstyle 2\times 2}  &  I_{\scriptscriptstyle 2\times 2} & \cdots & I_{\scriptscriptstyle 2\times 2}\end{array}\right]_{\scriptscriptstyle 2 \times 2(N-2)} 
 \end{split}
\end{equation*}
and
\begin{equation*}
\begin{split}
D_{\scriptscriptstyle 2 \times 2(N-2)} \left(B_{\scriptscriptstyle 2 \times 2}\right)^{-1}F_{\scriptscriptstyle 2 \times 2(N-2)}&=\left[ \begin{array}{c}
0_{\scriptscriptstyle 2\times 2}   \\
 0_{\scriptscriptstyle 2\times 2} 	\\
 \vdots \\
 -A_{\scriptscriptstyle 2\times 2} 
\end{array}\right]_{\scriptscriptstyle  2(N-2)\times 2} \left[\begin{array}{ccccc}
 I_{\scriptscriptstyle 2\times 2} &  I_{\scriptscriptstyle 2\times 2}  &  I_{\scriptscriptstyle 2\times 2} & \cdots & I_{\scriptscriptstyle 2\times 2}\end{array}\right]_{\scriptscriptstyle 2 \times 2(N-2)} \\
& = \left[ \begin{array}{cccc}
0_{\scriptscriptstyle 2\times 2}  & 0_{\scriptscriptstyle 2\times 2}  &\cdots & 0_{\scriptscriptstyle 2\times 2}   \\
0_{\scriptscriptstyle 2\times 2}  & 0_{\scriptscriptstyle 2\times 2}  &\cdots & 0_{\scriptscriptstyle 2\times 2}   \\
 \vdots &  \vdots &	 \cdots &  \vdots   \\
 0_{\scriptscriptstyle 2\times 2}  & 0_{\scriptscriptstyle 2\times 2}  &\cdots & 0_{\scriptscriptstyle 2\times 2}   \\
 -A_{\scriptscriptstyle 2\times 2}  &  -A_{\scriptscriptstyle 2\times 2}  &  \cdots &  -A_{\scriptscriptstyle 2\times 2} 
\end{array}\right]_{\scriptscriptstyle 2(N-2) \times 2(N-2)}.
 \end{split}
\end{equation*}
Therefore, we get
\begin{equation*}
C_{\scriptscriptstyle 2(N-2)\times 2(N-2)}  - D_{\scriptscriptstyle 2 \times 2(N-2)} \left(B_{\scriptscriptstyle 2 \times 2}\right)^{-1}F_{\scriptscriptstyle 2 \times 2(N-2)}=
\left[ \begin{array}{ccccc}
A_{\scriptscriptstyle 2\times 2} &-A_{\scriptscriptstyle 2\times 2}  & 0_{\scriptscriptstyle 2\times 2} & \cdots &0_{\scriptscriptstyle 2\times 2}  \\
 0_{\scriptscriptstyle 2\times 2} & A_{\scriptscriptstyle 2\times 2} &-A_{\scriptscriptstyle 2\times 2}  &  \cdots & 0_{\scriptscriptstyle 2\times 2} 	\\
 \vdots &  \vdots &   \vdots &   \vdots &   \vdots  \\
 A_{\scriptscriptstyle 2\times 2} & A_{\scriptscriptstyle 2\times 2} & A_{\scriptscriptstyle 2\times 2} & \cdots & 2A_{\scriptscriptstyle 2\times 2} 
\end{array}\right]_{\scriptscriptstyle 2(N-2) \times 2(N-2)}.
\end{equation*}
Then, $C_{\scriptscriptstyle 2(N-2)\times 2(N-2)}  - D_{\scriptscriptstyle 2 \times 2(N-2)} \left(B_{\scriptscriptstyle 2 \times 2}\right)^{-1}F_{\scriptscriptstyle 2 \times 2(N-2)}$  only depends of $A_{\scriptscriptstyle 2\times 2}$. Consequently, if $$\det \left( C_{\scriptscriptstyle 2(N-2)\times 2(N-2)}  - D_{\scriptscriptstyle 2 \times 2(N-2)} \left(B_{\scriptscriptstyle 2 \times 2}\right)^{-1}F_{\scriptscriptstyle 2 \times 2(N-2)}\right)=0,$$ we have that 
\begin{equation}\label{Ker}
\text{dim}\ \text{Ker} \left( C_{\scriptscriptstyle 2(N-2)\times 2(N-2)}  - D_{\scriptscriptstyle 2 \times 2(N-2)} \left(B_{\scriptscriptstyle 2 \times 2}\right)^{-1}F_{\scriptscriptstyle 2 \times 2(N-2)}\right)>0.
\end{equation}
\eqref{Ker}  implies that there exists a vector $$X_{\scriptscriptstyle 2(N-2)\times 1} \in \text{Ker} \left( C_{\scriptscriptstyle 2(N-2)\times 2(N-2)}  - D_{\scriptscriptstyle 2 \times 2(N-2)} \left(B_{\scriptscriptstyle 2 \times 2}\right)^{-1}F_{\scriptscriptstyle 2 \times 2(N-2)}\right)$$ such that 
\begin{equation}\label{new2}
X_{\scriptscriptstyle 2(N-2)\times 1}=(x_1,x_2, x_3, x_4 \cdots, x_{\scriptscriptstyle  2(N-3)}, x_{\scriptscriptstyle  2(N-2)})^T\neq 0_{\scriptscriptstyle 2(N-2)\times 1}
\end{equation}
and 
\begin{equation*}
\left( C_{\scriptscriptstyle 2(N-2)\times 2(N-2)}  - D_{\scriptscriptstyle 2 \times 2(N-2)} \left(B_{\scriptscriptstyle 2 \times 2}\right)^{-1}F_{\scriptscriptstyle 2 \times 2(N-2)}\right)\cdot X_{\scriptscriptstyle 2(N-2)\times 1} \\
=0_{\scriptscriptstyle 2(N-2) \times 1},
\end{equation*}
or equivalent,
\begin{equation}\label{new1}
\left[ \begin{array}{ccccc}
A_{\scriptscriptstyle 2\times 2} &-A_{\scriptscriptstyle 2\times 2}  & 0_{\scriptscriptstyle 2\times 2} & \cdots &0_{\scriptscriptstyle 2\times 2}  \\
 0_{\scriptscriptstyle 2\times 2} & A_{\scriptscriptstyle 2\times 2} &-A_{\scriptscriptstyle 2\times 2}  &  \cdots & 0_{\scriptscriptstyle 2\times 2} 	\\
 \vdots &  \vdots &   \vdots &   \vdots &   \vdots  \\
 A_{\scriptscriptstyle 2\times 2} & A_{\scriptscriptstyle 2\times 2} & A_{\scriptscriptstyle 2\times 2} & \cdots & 2A_{\scriptscriptstyle 2\times 2} 
\end{array}\right]_{\scriptscriptstyle 2(N-2) \times 2(N-2)}\left[ \begin{array}{ccccc}
x_{\scriptscriptstyle  1} \\
x_{\scriptscriptstyle  2}\\
x_{\scriptscriptstyle  3} \\
x_{\scriptscriptstyle 4} \\
\vdots \\
x_{\scriptscriptstyle  2(N-3)}\\
x_{ \scriptscriptstyle 2(N-2)}
\end{array}\right]_{\scriptscriptstyle 2(N-2) \times 1}=0_{\scriptscriptstyle 2(N-2) \times 1}.
\end{equation}
To finalize the proof of the Claim \ref{claim3}, denote 
\begin{equation*}
H^{\scriptscriptstyle 1}_{\scriptscriptstyle 2 \times 1}= \left[ \begin{array}{c}
x_{\scriptscriptstyle 1} \\
x_{\scriptscriptstyle 2}
\end{array}\right], \quad H^{\scriptscriptstyle 2}_{\scriptscriptstyle 2 \times 1}= \left[ \begin{array}{c}
x_{\scriptscriptstyle 3} \\
x_{\scriptscriptstyle 4}
\end{array}\right], \quad ..., \quad 
H^{\scriptscriptstyle  (N-2)}_{\scriptscriptstyle 2 \times 1}= \left[ \begin{array}{c}
x_(\scriptscriptstyle 2(N-3)) \\
x_(\scriptscriptstyle 2(N-2))
\end{array}\right].
\end{equation*}
Therefore, the product \eqref{new1} can be write in the form
\begin{equation*}
\left[ \begin{array}{ccccc}
A_{\scriptscriptstyle 2\times 2} &-A_{\scriptscriptstyle 2\times 2}  & 0_{\scriptscriptstyle 2\times 2} & \cdots &0_{\scriptscriptstyle 2\times 2}  \\
 0_{\scriptscriptstyle 2\times 2} & A_{\scriptscriptstyle 2\times 2} &-A_{\scriptscriptstyle 2\times 2}  &  \cdots & 0_{\scriptscriptstyle 2\times 2} 	\\
 \vdots &  \vdots &   \vdots &   \vdots &   \vdots  \\
 A_{\scriptscriptstyle 2\times 2} & A_{\scriptscriptstyle 2\times 2} & A_{\scriptscriptstyle 2\times 2} & \cdots & 2A_{\scriptscriptstyle 2\times 2} 
\end{array}\right]_{\scriptscriptstyle 2(N-2) \times 2(N-2)}\left[ \begin{array}{ccccc}
H^{\scriptscriptstyle  1}_{\scriptscriptstyle 2 \times 1} \\
H^{\scriptscriptstyle  2}_{\scriptscriptstyle 2 \times 1}\\
\vdots \\
H^{\scriptscriptstyle  (N-2)}_{\scriptscriptstyle 2 \times 1}
\end{array}\right]_{\scriptscriptstyle 2(N-2) \times 1}=0_{\scriptscriptstyle 2(N-2) \times 1}.
\end{equation*}
Thus, we have that 
\begin{equation*}
\left[ \begin{array}{ccccc}
A_{\scriptscriptstyle 2\times 2} \left( H^{\scriptscriptstyle  1}_{\scriptscriptstyle 2 \times 1} - H^{\scriptscriptstyle  2}_{\scriptscriptstyle 2 \times 1}\right)  \\
\\
A_{\scriptscriptstyle 2\times 2} \left( H^{\scriptscriptstyle  2}_{\scriptscriptstyle 2 \times 1} - H^{\scriptscriptstyle  3}_{\scriptscriptstyle 2 \times 1}\right)  \\
\\
\vdots \\
\\
A_{\scriptscriptstyle 2\times 2} \left( H^{\scriptscriptstyle  (N-1)}_{\scriptscriptstyle 2 \times 1} - H^{\scriptscriptstyle  (N-2)}_{\scriptscriptstyle 2 \times 1}\right)  \\
\\
A_{\scriptscriptstyle 2\times 2} \left( H^{\scriptscriptstyle  1}_{\scriptscriptstyle 2 \times 1} + H^{\scriptscriptstyle  2}_{\scriptscriptstyle 2 \times 1} +\cdots + 2H^{\scriptscriptstyle (N- 2)}_{\scriptscriptstyle 2 \times 1} \right)  
\end{array}\right]_{\scriptscriptstyle 2(N-2) \times 1}=0_{\scriptscriptstyle 2(N-2) \times 1}.
\end{equation*}
Now, let us now argue by contradiction.  If there exists $k \in \{ 1, 2,\cdots N-1\}$ such that $$H^{\scriptscriptstyle  k}_{\scriptscriptstyle 2 \times 1} -H^{\scriptscriptstyle  (k+1)}_{\scriptscriptstyle 2 \times 1} \neq 0_{\scriptscriptstyle 2\times 1},$$ we obtain
\begin{equation*}
A_{\scriptscriptstyle 2\times 2} \left( H^{\scriptscriptstyle  k}_{\scriptscriptstyle 2 \times 1} -H^{\scriptscriptstyle  (k+1)}_{\scriptscriptstyle 2 \times 1}\right) = 0_{\scriptscriptstyle 2\times 1},
\end{equation*}
which implies that dim $	\text{ker} A_{\scriptscriptstyle 2\times 2} >0$, it means that $\det A_{\scriptscriptstyle 2\times 2} =0$. However, from the case $N=2$, we known that 
\begin{equation*}
\det A_{\scriptscriptstyle 2\times 2} = ag -fn \neq 0,
\end{equation*} 
and hence we obtain a contradiction. On the other hand, suppose that 
\begin{equation}\label{new3}
H^{\scriptscriptstyle  j}_{\scriptscriptstyle 2 \times 1} -H^{\scriptscriptstyle  (j+1)}_{\scriptscriptstyle 2 \times 1} = 0_{\scriptscriptstyle 2\times 1}, \quad \forall j=1,2,..., N-2.
\end{equation}
Thus, from 	\eqref{new2} and \eqref{new3}, we deduce that $H^{\scriptscriptstyle  j}_{\scriptscriptstyle 2 \times 1}\neq 0_{\scriptscriptstyle 2\times 1}$ for some $j \in \{1,2,..., N-2\}$ and 
\begin{equation*}
A_{\scriptscriptstyle 2\times 2} \left( H^{\scriptscriptstyle  1}_{\scriptscriptstyle 2 \times 1} + H^{\scriptscriptstyle  2}_{\scriptscriptstyle 2 \times 1} +\cdots + 2H^{\scriptscriptstyle (N- 2)}_{\scriptscriptstyle 2 \times 1} \right)  = (N-1)A_{\scriptscriptstyle 2\times 2} H^{\scriptscriptstyle  j}_{\scriptscriptstyle 2 \times 1} =  0_{\scriptscriptstyle 2(N-2) \times 1}.
\end{equation*}
Which is again a contradiction, by using the case $N=2$. Hence, in the two cases, we only have that 
$$\det \left( C_{\scriptscriptstyle 2(N-2)\times 2(N-2)}  - D_{\scriptscriptstyle 2 \times 2(N-2)} \left(B_{\scriptscriptstyle 2 \times 2}\right)^{-1}F_{\scriptscriptstyle 2 \times 2(N-2)}\right)\neq 0,$$
it implies that  $\det \mathbf{M'}\neq 0$. Consequently, the determinant of $\mathbf{M}$ is nonzero, implying that the matrix $\mathbf{M}$ is invertible and the Claim \ref{claim3} follows. Thus, Lemma \ref{tec} is proved. \qed

\appendix

\section{Vertex conditions types $\mathcal{B}$ and $\mathcal{C}$}\label{apenB}
In this appendix, we will outline how to prove that matrices associated with vertex conditions \eqref{boundary2} (type $\mathcal{B}$) and \eqref{boundary3} (type $\mathcal{C}$)  are invertible. We will  consider the vertex conditions 
\begin{equation*}
\text{Type $\mathcal{B}$: }\begin{cases}
\partial_x^k u_1(t,0)=\partial_x^k u_2(t,0)=\cdots=	\partial_x^k u_N(t,0),\ k=2,3 & t\in(0,T),\\
\sum_{j=1}^N\partial_x^ku_j(t,0)=0,\ k=0,1\ & t\in(0,T),
\end{cases}
\end{equation*}
and
\begin{equation*}
\text{Type $\mathcal{C}$: }\begin{cases}
\partial_x^k u_1(t,0)=\partial_x^k u_2(t,0)=\cdots=\partial_x^k u_N(t,0),\ k=0,3 & t\in(0,T),\\
\sum_{j=1}^N\partial_x^ku_j(t,0)=0,\ k=1,2\ & t\in(0,T),
\end{cases}
\end{equation*}
which may be expressed in matrices form as follows
\begin{align*}
\mathbf{M_{\mathcal{B}}}\,\, \boldsymbol{\gamma_{\mathcal{B}}}=\mathbf{F_{\mathcal{B}}}, 	\quad \text{and} 	\quad \mathbf{M_{\mathcal{C}}}\,\,\boldsymbol{\gamma_{\mathcal{C}}}=\mathbf{F_{\mathcal{C}}}, 
\end{align*}
respectively. Here  $\mathbf{F_{\mathcal{B}}}$ and $\mathbf{F_{\mathcal{C}}}$ are the function vectors defined by
\begin{align*}
\mathbf{F_{\mathcal{B}}}:=-\left[\begin{array}{c}
\sum_{j=1}^N F_j(t,0)\\
\sum_{j=1}^N \partial_x\mathcal{I}_{\frac14}F_j(t,0)\\
\\
\partial_x^2\mathcal{I}_{\frac12}F_1(t,0)-\partial_x^2\mathcal{I}_{\frac12}F_2(t,0)\\
\vdots \\
\partial_x^2\mathcal{I}_{\frac12}F_{N-1}(t,0)-\partial_x^2\mathcal{I}_{\frac12}F_{N}(t,0)\\
\\
\partial_x^3\mathcal{I}_{\frac34}F_1(t,0)-\partial_x^3\mathcal{I}_{\frac34}F_2(t,0)\\
\vdots \\
\partial_x^3\mathcal{I}_{\frac34}F_{N-1}(t,0)-\partial_x^3\mathcal{I}_{\frac34}F_{N}(t,0)\\
\end{array}\right]_{\scriptscriptstyle 2N\times 1} \quad 
\mathbf{F_{\mathcal{C}}}:=-\left[\begin{array}{c}
F_1(t,0)-F_2(t,0)\\
\vdots \\
F_{N-1}(t,0)-F_N(t,0)\\
\\
\sum_{j=1}^N \partial_x\mathcal{I}_{\frac14}F_j(t,0)\\
\sum_{j=1}^N \partial_x^2\mathcal{I}_{\frac12}F_j(t,0)\\
\\
\partial_x^3\mathcal{I}_{\frac34}F_1(t,0)-\partial_x^3\mathcal{I}_{\frac34}F_2(t,0)\\
\vdots \\
\partial_x^3\mathcal{I}_{\frac34}F_{N-1}(t,0)-\partial_x^3\mathcal{I}_{\frac34}F_{N}(t,0)\\
\end{array}\right]_{\scriptscriptstyle 2N\times 1},
\end{align*}
$\boldsymbol{\gamma_{\mathcal{B}}}$ and $\boldsymbol{\gamma_{\mathcal{C}}}$  are the  matrices column given by vectors $(\gamma_{11}^{\mathcal{B}},\gamma_{12}^{\mathcal{B}},\cdots,\gamma_{N1}^{\mathcal{B}},\gamma_{N2}^{\mathcal{B}})$ and $(\gamma_{11}^{\mathcal{C}},\gamma_{12}^{\mathcal{C}},\cdots,\gamma_{N1}^{\mathcal{C}},\gamma_{N2}^{\mathcal{C}})$, respectively. 

Note that choosing $\lambda_{11}=\lambda_{21} =\cdots =\lambda_{N1}$ and $\lambda_{12}=\lambda_{22}=\cdots=\lambda_{N2}$ and arguing as in the Section \ref{sec5}, the determinants of the matrices 
\begin{align*}
\mathbf{M}_{\mathcal{B}}=\mathbf{M}_{\mathcal{B}}(-1/2,1/4,\cdots, -1/2,1/4)\quad \text{ and} 	\quad \mathbf{M}_{\mathcal{C}}=\mathbf{M}_{\mathcal{C}}(-1/2,1/4,\cdots, -1/2,1/4)
\end{align*}
 are given by 
\begin{align*}
|\mathbf{ M_{\mathcal{B}}}|&= \frac{M^{2N}e^{\frac{ (12+N)i\pi}{8}}}{8^{2N}}  \prod_{i=1,2, j=1,\cdots,N}  \left( e^{-\frac{3i\pi\lambda_{ji}}{8}}+e^{\frac{5i\pi\lambda_{ji}}{8}}\right)|\mathbf{ M'_{\mathcal{B}}}|,
\end{align*}
and 
\begin{align*}
|\mathbf{ M_{\mathcal{C}}}|&= \frac{M^{2N}e^{\frac{ (12+N)i\pi}{8}}}{8^{2N}}  \prod_{i=1,2, j=1,\cdots,N}  \left( e^{-\frac{3i\pi\lambda_{ji}}{8}}+e^{\frac{5i\pi\lambda_{ji}}{8}}\right)|\mathbf{ M'_{\mathcal{C}}}|,
\end{align*}
where
\begin{equation*}
\left| \mathbf{M'_{\mathcal{B}}} \right|=\left|\begin{array}{rrrrrrrrrrr}
	a& 	n& a&  n & a & n & \cdots &a &n &a&n\\
	f& 	g& f&g & f & g & \cdots &f &g &f&g\\
	\\
	c& 	e& -c&-e& 0 & 0 & \cdots &0 &0 &0&0\\
	0&0&	c& 	e& -c&-e &  \cdots &0 &0 &0&0\\
	\vdots & 	\vdots & 	\vdots & 	\vdots & 	\vdots & 	\vdots & 	\vdots & 	\vdots & 	\vdots & 	\vdots & 	\vdots \\
	0&0&0&0&0&0&\cdots&	c& e& -c&-e \\
	\\
	d& 	m& -d&-m& 0 & 0 & \cdots &0 &0 &0&0\\
	0&0&	d& 	m& -d&-m &  \cdots &0 &0 &0&0\\
	\vdots & 	\vdots & 	\vdots & 	\vdots & 	\vdots & 	\vdots & 	\vdots & 	\vdots & 	\vdots & 	\vdots & 	\vdots \\
	0&0&0&0&0&0&\cdots&	d& m& -d&-m \\
	\end{array}\right|_{\scriptscriptstyle 2N\times 2N},
\end{equation*}
and 
\begin{equation*}
\left| \mathbf{M'_{\mathcal{C}}} \right|=\left|\begin{array}{rrrrrrrrrrr}
	a& 	n& -a&-n& 0 & 0 & \cdots &0 &0 &0&0\\
	0&0&	a& 	n& -a&-n&  \cdots &0 &0 &0&0\\
	\vdots & 	\vdots & 	\vdots & 	\vdots & 	\vdots & 	\vdots & 	\vdots & 	\vdots & 	\vdots & 	\vdots & 	\vdots \\
	0&0&0&0&0&0&\cdots&	a& n& -a&-n \\
	\\
	f& 	g& f&g & f & g & \cdots &f &g &f&g\\
	c& 	e& c&  e & c & e & \cdots &c &e &c&e\\
	\\
	d& 	m& -d&-m& 0 & 0 & \cdots &0 &0 &0&0\\
	0&0&	d& 	m& -d&-m &  \cdots &0 &0 &0&0\\
	\vdots & 	\vdots & 	\vdots & 	\vdots & 	\vdots & 	\vdots & 	\vdots & 	\vdots & 	\vdots & 	\vdots & 	\vdots \\
	0&0&0&0&0&0&\cdots&	d& m& -d&-m \\
	\end{array}\right|_{\scriptscriptstyle 2N\times 2N}.
\end{equation*}

As in the case vertex type $\mathcal{A}$, we need to study the determinants of the matrices $\mathbf{M'_{\mathcal{B}}}$ and $\mathbf{M'_{\mathcal{C}}}$. In order to see the determinant of  $\mathbf{M'_{\mathcal{B}}}$ and $\mathbf{M'_{\mathcal{C}}}$ are nonzero, we use the determinant properties together with \eqref{blockdeter} to observe that
\begin{equation*}
\det \mathbf{M'_{\mathcal{B}}} =\det \mathbf{M''_{\mathcal{B}}}  \det \left[ \begin{array}{c c}
a & n \\
f & g
\end{array}\right] \quad \text{and} \quad \det \mathbf{M'_{\mathcal{C}}} =\det \mathbf{M''_{\mathcal{C}}}  \det \left[ \begin{array}{c c}
c & e \\
f & g
\end{array}\right].
\end{equation*}
These two matrices, namely,  $\mathbf{M'_{\mathcal{B}}} $ and $ \mathbf{M'_{\mathcal{C}}} $, have the following two following properties:
\begin{enumerate}
\item[(i)] If $det \left[ \begin{array}{c c}
c & e \\
d & m
\end{array}\right]=cm-de\neq 0$, then $ \det \mathbf{M''_{\mathcal{B}}}\neq 0$.

\item[(ii)] If $det \left[ \begin{array}{c c}
a & n \\
d& m
\end{array}\right]=am-dn\neq 0$, then $ \det \mathbf{M''_{\mathcal{C}}}\neq 0$.
\end{enumerate}

\begin{claim}\label{claim4}
The relations
$$ag -fn \neq, \quad cg - fe \neq 0$$
and
$$cm-de\neq 0, \quad am-dn\neq 0$$
%\begin{align*}
%\begin{cases}
%ag -fn \neq 0,  & cg - fe \neq 0 \\
%cm-de\neq 0, &  am-dn\neq 0
%\end{cases}
%\end{align*}
are valid.
\end{claim}

In fact, choosing
\begin{align*}
\lambda_{11}=\lambda_{21}=\cdots=\lambda_{N1}=-\frac12 \quad \text{and} \quad \lambda_{11}=\lambda_{21}=\cdots=\lambda_{N1}=-\frac14,
\end{align*}
together with  \eqref{det1}, \eqref{new6} and \eqref{new5}, we already now that $ag -fn$ and $cm-de$ are nonzero. Finally,  easy calculations yield that
\begin{align*}
cg-fe &= \left( \frac{1}{\sin\left(\frac78\pi\right)} \right)\left(-\frac{\tan\left(\frac{\pi}{8}\right)}{	\cos\left(\frac{\pi	}{16}\right)	}\right) - \left( \frac{1}{\sin\left(\frac{11}{16}\pi\right)} \right)\left(-\frac{\tan\left(-\frac{\pi}{4}\right)}{	\cos\left(-\frac{\pi	}{8}\right)	}\right) \sim -2.4053\neq 0,  \\
\\
am-dn&=\left( \frac{1}{\sin\left(\frac38\pi\right)} \right)\left(-\frac{\tan\left(\frac{\pi}{8}\right)}{	\sin\left(\frac{\pi	}{16}\right)	}\right) - \left( \frac{1}{\sin\left(\frac{5}{16}\pi\right)} \right)\left(-\frac{\tan\left(-\frac{\pi}{4}\right)}{	\sin\left(-\frac{\pi	}{8}\right)	}\right)  \sim 0.8446\neq 0,
\end{align*}
showing the Claim \ref{claim4}, and thus the matrices $\mathbf{M_{\mathcal{B}}}$ and $\mathbf{M_{\mathcal{C}}}$ are invertible. 

\subsection{Proof of Theorem \ref{theorem1}: Vertex conditions type $\mathcal{B}$ and $\mathcal{C}$}
The analysis developed above give us the following representations  for $\boldsymbol{\gamma_{\mathcal{B}}}$ and $\boldsymbol{\gamma_{\mathcal{C}}}$
\begin{align*}
\boldsymbol{\gamma_{\mathcal{B}}}=\mathbf{M_{\mathcal{B}}}^{-1}\,\, \mathbf{F_{\mathcal{B}}} 	\quad \text{and} 	\quad \boldsymbol{\gamma_{\mathcal{C}}}=\mathbf{M_{\mathcal{C}}}^{-1}\,\,\mathbf{F_{\mathcal{C}}},
\end{align*}
respectively. Therefore, the solution $u_j(t,x)$ of the Cauchy problem \eqref{grapha} with vertex conditions type $\mathcal{B}$ and $\mathcal{C}$ can be express in a integral forms
\begin{equation}\label{form2}
u_j^{\mathcal{B}}(t,x)=\mathcal{L}^{-\frac12}\gamma_{j1}^{\mathcal{B}}(t,x)+\mathcal{L}^{\frac14}\gamma_{j2}^{\mathcal{B}}(t,x)+F_j(t,x),
\quad j=1,2, ..., N
\end{equation}
and
\begin{equation}\label{form3}
u_j^{\mathcal{C}}(t,x)=\mathcal{L}^{-\frac12}\gamma_{j1}^{\mathcal{C}}(t,x)+\mathcal{L}^{\frac14}\gamma_{j2}^{\mathcal{C}}(t,x)+F_j(t,x),
\quad j=1,2, ..., N.
\end{equation} 

Finally, in order to establish Theorem \ref{theorem1} with boundary conditions type $\mathcal{B}$ and $\mathcal{C}$, we closely follow the same steps of subsection \ref{truncated}, \ref{welldefined} and \ref{contraction}, for use the Fourier restriction method to define a truncated version for \eqref{form2} and \eqref{form3}, proving thus that $\mathcal{L}^{-\frac12}\gamma_{j1}^{\mathcal{B}}$, $\mathcal{L}^{-\frac12}\gamma_{j1}^{\mathcal{C}}$, $\mathcal{L}^{\frac14}\gamma_{j2}^{\mathcal{B}}$ and $\mathcal{L}^{-\frac12}\gamma_{j1}^{\mathcal{C}}$  for $j=1,2,..., N$, are well-defined. With this in hand, a contraction mapping argument gives us the result desired. \qed

\subsection*{Acknowledgments} %The authors wish to thank the referees for his/her valuable comments which improved this paper. 
Capistrano-Filho was supported by CNPq 306475/2017-0, 408181/2018-4, CAPES-PRINT 88881.311964/2018-01 and Qualis A - Propesq (UFPE).  Gallego was partially supported by Facultad de Ciencias Exactas y Naturales, Nacional de Colombia Sede Manizales, under the project 45511 and CAPES-PRINT under grant number 88881.311964/2018-01. 

This work was carried out during some visits of the authors to the Federal University of Pernambuco, Federal University of Alagoas and Universidad Nacional de Colombia - Sede Manizales. The authors would like to thank the Universities for its hospitality. 

\end{document}